\documentclass[12pt]{amsart}
\usepackage{amsfonts,amssymb,amscd,amsmath,enumerate,color, xypic, fancybox, graphics}
\def\frk{\mathfrak}               

\def\Phi{{\frk N}}
\def\opn#1#2{\def#1{\operatorname{#2}}} 
\opn\chara{char} 
\opn\length{\ell} 
\opn\pd{pd} 
\opn\rk{rk}
\opn\projdim{proj\,dim} 
\opn\injdim{inj\,dim} 
\opn\rank{rank}
\opn\depth{depth} 
\opn\grade{grade} 
\opn\height{height}
\opn\embdim{emb\,dim} 
\opn\codim{codim}

\opn\Tr{Tr} 
\opn\bigrank{big\,rank}
\opn\superheight{superheight}
\opn\lcm{lcm}
\opn\trdeg{tr\,deg}
\opn\reg{reg} 
\opn\lreg{lreg} 
\opn\ini{in} 
\opn\lpd{lpd}
\opn\size{size}
\opn\mult{mult}
\opn\dist{dist}
\opn\cone{cone}
\opn\lex{lex}
\opn\rev{rev}
\opn\div{div} \opn\Div{Div} \opn\cl{cl} \opn\Cl{Cl}
\opn\Spec{Spec} \opn\Supp{Supp} \opn\supp{supp} \opn\Sing{Sing}
\opn\Ass{Ass} \opn\Min{Min}
\opn\Ann{Ann} \opn\Rad{Rad} \opn\Soc{Soc}
\opn\Syz{Syz} \opn\Im{Im} \opn\Ker{Ker} \opn\Coker{Coker}
\opn\Am{Am} \opn\Hom{Hom} \opn\Tor{Tor} \opn\Ext{Ext}
\opn\End{End} \opn\Aut{Aut} \opn\id{id} \opn\ini{in}

\opn\nat{nat}
\opn\pff{pf}
\opn\Pf{Pf} \opn\GL{GL} \opn\SL{SL} \opn\mod{mod} \opn\ord{ord}
\opn\Gin{Gin}
\opn\Hilb{Hilb}\opn\adeg{adeg}\opn\std{std}\opn\ip{infpt}
\opn\Pol{Pol}
\opn\sat{sat}
\opn\Var{Var}
\opn\Gen{Gen}

\opn\aff{aff} \opn\con{conv} \opn\relint{relint} \opn\st{st}
\opn\lk{lk} \opn\cn{cn} \opn\core{core} \opn\vol{vol}
\opn\link{link} \opn\star{star}
\opn\gr{gr}

\def\pot#1#2{#1[\kern-0.28ex[#2]\kern-0.28ex]}

\opn\dirlim{\underrightarrow{\lim}}
\opn\inivlim{\underleftarrow{\lim}}
\let\to=\rightarrow

\def\Implies{\ifmmode\Longrightarrow \else
        \unskip${}\Longrightarrow{}$\ignorespaces\fi}
\def\implies{\ifmmode\Rightarrow \else
        \unskip${}\Rightarrow{}$\ignorespaces\fi}
\def\iff{\ifmmode\Longleftrightarrow \else
        \unskip${}\Longleftrightarrow{}$\ignorespaces\fi}

\let\:=\colon
\newtheorem{Theorem}{Theorem}[section]
\newtheorem{Lemma}[Theorem]{Lemma}
\newtheorem{Corollary}[Theorem]{Corollary}
\newtheorem{Proposition}[Theorem]{Proposition}
\newtheorem{Remark}[Theorem]{Remark}

\newtheorem{Example}[Theorem]{Example}

\newtheorem{Definition}[Theorem]{Definition}

\newtheorem{Conjecture}[Theorem]{Conjecture}
\newtheorem{Question}[Theorem]{Question}
\let\epsilon\varepsilon
\let\phi=\varphi
\let\kappa=\varkappa
\def\qed{\ifhmode\textqed\fi
      \ifmmode\ifinner\quad\qedsymbol\else\dispqed\fi\fi}
\def\textqed{\unskip\nobreak\penalty50
       \hskip2em\hbox{}\nobreak\hfil\qedsymbol
       \parfillskip=0pt \finalhyphendemerits=0}
\def\dispqed{\rlap{\qquad\qedsymbol}}

\opn\dis{dis}
\opn\height{height}
\opn\dist{dist}
\def\pnt{{\raise0.5mm\hbox{\large\bf.}}}

\opn\Lex{Lex}
\opn\conv{conv}

\begin{document}
\title{Weakly closed graphs and F-purity of \\ binomial edge ideals}
\author[K.~Matsuda]{Kazunori Matsuda}
\address[Kazunori Matsuda]{Department of Pure and Applied Mathematics,
	Graduate School of Information Science and Technology,
	Osaka University,
	Suita, Osaka 565-0871, Japan}
\email{kaz-matsuda@ist.osaka-u.ac.jp}
\begin{abstract}
Herzog-Hibi-Hreind\'{o}ttir-Kahle-Rauh introduced the class of closed graph and 
they proved that the binomial edge ideal $J(G)$ of a graph $G$ has quadratic Gr\"{o}bner bases 
if $G$ is closed. 
In this paper, we introduce the class of weakly closed graph as a generalization of the closed graph 
and prove that the quotient ring $S/J(G)$ is $F$-pure if $G$ is weakly closed. 
This fact is a generalization of Ohtani's theorem. 
\end{abstract}
\maketitle

\section*{Introduction}
Let $G$ be a simple (i.e. $G$ has no loops and multiple edges) graph on the vertex set 
$V(G) = [n] = \{1, 2, \ldots, n\}$ with the edge set $E(G)$. 
Let $S = K[x_1, \ldots, x_n, y_1, \ldots, y_n]$ be the polynomial ring over a field $K$. 

The following binomial ideal $J_{G}$ of $S$
\[
J_{G} := ([i, j] = x_{i}y_{j} - x_{j}y_{i} \mid \{i, j\} \in E(G) \ \text{with} \  i < j)
\] 
is called the {\em binomial edge ideal} of $G$. 
This was introduced by Herzog-Hibi-Hreind\'{o}ttir-Kahle-Rauh \cite{origin} and 
Ohtani \cite{O1} independently. 
The binomial edge ideal has been studied multilaterally by many researchers, e.g:
\begin{itemize}
	\item Cohen-Macaulayness (\cite{BMS}, \cite{CMBEI}, \cite{origin}, \cite{KS}, \cite{RaRi}, \cite{Rinaldo}, \cite{cactus}, \cite{Z}),
	\item betti number and regularity (\cite{B}, \cite{block}, \cite{deAH}, \cite{Dokuyucu}, \cite{EZ}, \cite{pure}, \cite{regularity}, \cite{MM}, \cite{SK}, \cite{reg3}, \cite{SZ}, \cite{ZZ}),
	\item Koszulness (\cite{KoszulBEI}, \cite{KF}, \cite{dual}), 
	\item Gr\"{o}bner bases (\cite{BBS}, \cite{CrRi1}, \cite{origin}, \cite{O1}), 
	\item generalizations of the binomial edge ideal (\cite{CDG}, \cite{EHHQ}, \cite{R}, \cite{pair}), 
\end{itemize} 
etc. 

The notion of closed graph, introduced by Herzog et al. (\cite[p.319]{origin}), 
appears in many of studies of binomial edge ideals. 
A simple graph $G$ is {\em closed with respect to the given labeling of the vertices} 
if the condition is satisfied: for all $\{i, j\}, \{k, l\} \in E(G)$ with $i < j$ and $k < l$, 
one has $\{j, l\} \in E(G)$ if $i = k$, and $\{i, k\} \in E(G)$ if $j = l$. 
A simple  graph $G$ is {\em closed} if there exists a labeling such that it is closed. 
Crupi-Rinaldo proved that $J_{G}$ has a quadratic Gr\"{o}bner basis if and only if 
$G$ is closed (\cite[Theorem 3.4]{CrRi1}). 
In addition, graph-theoretic properties of closed graphs are also studied 
(e.g. \cite{closed1}, \cite{closed2}, \cite{Cr}, \cite{CrRi2})
and Sharifan-Javanbakht introduced the notion of $m$-closed graph 
as a generalization of closed graph (\cite[Definition 5]{mclosed}). 

In this paper, we introduce the notion of {\em weakly closed graph} 
as a generalization of closed graph, and we prove that the quotient ring 
$S/J_{G}$ is $F$-pure if $G$ is weakly closed. 

Let us explain the organization of this paper. 

In Section 1, we define the weakly closed graph and discuss the difference between 
the closed graph and the weakly closed graph. 
We also state some classes of graphs related to the weakly closed graph. 

In Section 2, we state the definition of $F$-purity and 
prove that $S/J_{G}$ is $F$-pure if $G$ is weakly closed. 

In Section 3, we classify all connected weakly closed graphs with $5$ and $6$ vertices 
by using Harary's classification of connected graph \cite{Ha}. 

About terminologies for the graph theory, see \cite{Diestel}. 

\section{Weakly closed graph}

Let $G$ be a simple graph on the vertex set 
$V(G) = [n]$ with the edge set $E(G)$. 

\begin{Definition}
$G$ is said to be {\em weakly closed} if there exists a labeling of the vertices 
which satisfies the following condition: 
for all integers $1 \le i < j < k \le n$, if $\{i, k\} \in E(G)$ then $\{i, j\} \in E(G)$ or $\{j ,k\} \in E(G)$. 
\end{Definition}

\begin{Example}
The following graph $G$ is weakly closed:

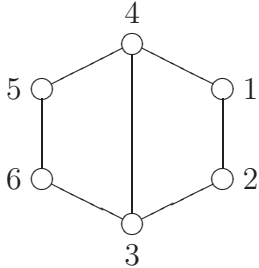
\begin{xy}
	\ar@{} (60, -20) *++!R{6} *\cir<4pt>{} = "A";
	\ar@{} "A";(84, -20) *++!L{2} *\cir<4pt>{} = "B";
	\ar@{-} "A";(72, -26) *++!U{3} *\cir<4pt>{} = "C";
	\ar@{-} "A";(60, -8) *++!R{5} *\cir<4pt>{} = "D";
	\ar@{-} "D";(72, -2) *++!D{4} *\cir<4pt>{} = "E";
	\ar@{} "D";(84, -8) *++!L{1} *\cir<4pt>{} = "F";
	\ar@{-} "B";"C";
	\ar@{-} "B";"F";
	\ar@{-} "E";"F";
	\ar@{-} "C";"E";
\end{xy}
\end{Example}

From the following proposition, we can see that the weakly closed graph is a 
generalization of the closed graph. 

\begin{Proposition}$($See \cite[Proposition 4.8]{CrRi2}$)$
Let $G$ be a graph. 
Then the following conditions are equivalent:
	\begin{enumerate}[$(1)$]
		\item $G$ is closed. 
		\item There exists a labeling of the vertices 
which satisfies the following condition: 
for all integers $1 \le i < j < k \le n$, if $\{i, k\} \in E(G)$ then $\{i, j\} \in E(G)$ and $\{j ,k\} \in E(G)$. 
	\end{enumerate}
\end{Proposition}
\begin{proof}
$(1) \Rightarrow (2)$: Take a closed labeling and assume that $\{i, k\} \in E(G)$.
Since $G$ is closed, one has $\{i, i + 1\}, \{i + 1, i + 2\}, \ldots, \{k - 1, k\} \in E(G)$ 
by \cite[Proposition 1.4]{origin}.  
Then we have that $\{i, k - 1\}, \{i, k - 2\}, \ldots, \{i, i + 2\} \in E(G)$ by the definition of closedness. 
Similarly, we also have that $\{i + 1, k\}, \{i + 2, k\}, \ldots, \{k - 2, k\} \in E(G)$. 

$(2) \Rightarrow (1)$: Assume that $i < j < k$. 
If $\{i, j\}, \{i, k\} \in E(G)$, then $\{j, k\} \in E(G)$ by assumption. 
Similarly, if $\{i, k\}, \{j, k\} \in E(G)$, then $\{i, j\} \in E(G)$. 
Therefore $G$ is closed. 
\end{proof}

Here we give some remarks about the difference between closed graphs and weakly closed graphs. 

\begin{Proposition}
Let $G$ be a connected graph. 
Then
	\begin{enumerate}[$(1)$]
		\item \cite[Proposition 1.4]{origin}  If $G$ is closed, then $\{i, i + 1\} \in E(G)$ 
		for {\em all} $1 \le i \le n - 1$. 
		\item If $G$ is weakly closed, then $\{i, i + 1\} \in E(G)$ 
		for {\em some} $1 \le i \le n - 1$.
	\end{enumerate}
\end{Proposition}
\begin{proof}
(2) : Assume that $G$ is weakly closed and $\{i, i + 1\} \not\in E(G)$ for all $1 \le i \le n - 1$. 
Since $\{1, 2\}, \{2, 3\} \not\in E(G)$, we have $\{1, 3\} \not\in E(G)$. 
Hence $\{1, 4\} \not\in E(G)$ from $\{3, 4\} \not\in E(G)$. 
By repeating this argument, we have $\{1, 2\}, \{1, 3\}, \ldots, \{1, n\} \not\in E(G)$, 
but this is a contradiction since $G$ is connected. 
Therefore we have the desired conclusion. 
\end{proof}

We call the following graphs {\em claw} and {\em bigclaw} respectively. \\

\begin{xy}
	\ar@{} (50, -8) *++!R{} *\cir<4pt>{} = "H";
	\ar@{-} "H";(50, 8) *++!R{} *\cir<4pt>{} = "I";
	\ar@{-} "H";(34, -8) *++!R{} *\cir<4pt>{} = "J";
	\ar@{-} "H";(66, -8) *++!R{} *\cir<4pt>{} = "K";
	\ar@{} "H";(50, -14) *++!U{{\rm claw}} = "L";	
	\ar@{} (102, -8) *++!R{} *\cir<4pt>{} = "A";
	\ar@{-} "A";(102, 0) *++!R{} *\cir<4pt>{} = "B";
	\ar@{-} "B";(102, 8) *++!R{} *\cir<4pt>{} = "C";
	\ar@{-} "A";(94, -8) *++!R{} *\cir<4pt>{} = "D";
	\ar@{-} "D";(86, -8) *++!R{} *\cir<4pt>{} = "E";
	\ar@{-} "A";(110, -8) *++!R{} *\cir<4pt>{} = "F";
	\ar@{-} "F";(118, -8) *++!R{} *\cir<4pt>{} = "G";
	\ar@{} "A";(102, -14) *++!U{{\rm bigclaw}} = "N";						
\end{xy}

\begin{Proposition}\label{claw}
Let $G$ be a graph. 
Then
	\begin{enumerate}[$(1)$]
		\item \cite[Proposition 1.2]{origin} If $G$ is closed, then $G$ is claw-free and chordal. 
		\item If $G$ is weakly closed, then $G$ is bigclaw-free and every cycle of length $5$ 
		or more in $G$ has a chord. 
	\end{enumerate}
\end{Proposition}
\begin{proof}
(2) : It is enough to show that the bigclaw and the chordless cycle $C_n$ of length $n \ge 5$ are not 
weakly closed. 

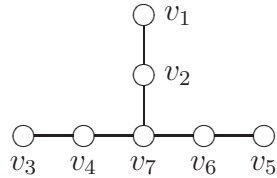
\vspace{3mm}

\begin{xy}
	\ar@{} (75, -8) *++!U{v_{7}} *\cir<4pt>{} = "A";
	\ar@{-} "A";(75, 0) *++!L{v_{2}} *\cir<4pt>{} = "B";
	\ar@{-} "B";(75, 8) *++!L{v_{1}} *\cir<4pt>{} = "C";
	\ar@{-} "A";(67, -8) *++!U{v_{4}} *\cir<4pt>{} = "D";
	\ar@{-} "D";(59, -8) *++!U{v_{3}} *\cir<4pt>{} = "E";
	\ar@{-} "A";(83, -8) *++!U{v_{6}} *\cir<4pt>{} = "F";
	\ar@{-} "F";(91, -8) *++!U{v_{5}} *\cir<4pt>{} = "G";						
\end{xy}

\vspace{3mm}

First, we show that the bigclaw is not weakly closed. 
Assume that the bigclaw $G$ is weakly closed under the above labeling 
and $V(G) = \{v_{1}, \ldots, v_{7}\} = \{1, \ldots, 7\}$. 
If $1 \in \{v_{1}, v_{2}\}$, then $\{v_{1}, v_{2}\} = \{1, 2\}$ or $\{v_{1}, v_{2}\} = \{1, 3\}$. 

\vspace{2mm}

Case 1: Assume that $\{v_{1}, v_{2}\} = \{1, 2\}$. 
Then $3 \le v_{7} \le 5$. 
	\begin{itemize}
		\item $v_{7} = 3 \Rightarrow 7 \not\in \{v_{4}, v_{6}\}$. 
		We may assume that $v_{3} = 7$.  Then $v_{4} = 6$, but this is a contradiction 
		since $3 < v_{5} < 6$, $\{3, 6\} \in E(G)$ but $\{3, v_{5}\}, \{v_{5}, 6\} \not\in E(G)$. 
		\item $v_{7} = 4 \Rightarrow 3 \in \{v_{4}, v_{6}\}$. 
		We may assume that $v_{4} = 3$. Then $v_{3} = 5$, but this is a contradiction 
		since $4 < 5 < v_{6}$, $\{4, v_{6}\} \in E(G)$ but $\{4, 5\}, \{5, v_{6}\} \not\in E(G)$. 
		\item $v_{7} = 5 \Rightarrow 3, 4 \in \{v_{4}, v_{6}\}$. 
		We may assume that $v_{4} = 3$ and $v_{6} = 4$. 
		Then $v_{3} \ge 6$, but this is a contradiction since $3 < 4 < v_{3}$, 
		$\{3, v_{3}\} \in E(G)$ but $\{3, 4\}, \{4, v_{3}\} \not\in E(G)$. 
	\end{itemize}
	
Case 2 : Assume that $\{v_{1}, v_{2}\} = \{1, 3\}$. Then $v_{7} = 2$, hence 
$7 \not\in \{v_{4}, v_{6}\}$. 
We may assume that $v_{3} = 7$. Then $v_{4} = 6$, but this is a contradiction 
since $2 < v_{5} < 6$, $\{2, 6\} \in E(G)$ but $\{2, v_{5}\}, \{v_{5}, 6\} \not\in E(G)$.   

From the above argument, we have $1 \not\in \{v_{1}, v_{2}\}$. 
Similarly, we also have $1 \not\in \{v_{3}, v_{4}, v_{5}, v_{6}\}$. 
Hence $v_{7} = 1$. Then $7 \not\in \{v_{2}, v_{4}, v_{6}\}$. 
We may assume that $v_{1} = 7$. 
Then $v_{2} = 6$, but this is a contradiction since $1 < v_{3} < 6$, $\{1, 6\} \in E(G)$ 
but $\{1, v_{3}\}, \{v_{3}, 6\} \not\in E(G)$. 
Therefore the bigclaw is not weakly closed. 

Next, we show that the chordless cycle $C_{n}$ of length $n \ge 5$ is not weakly closed. 
Assume that $C_{n}$ is weakly closed and $V(C_{n}) = \{v_{1}, \ldots, v_{n}\} = \{1, \ldots, n\}$, 
$E(C_{n}) = \{ \{v_{1}, v_{2}\}, \{v_{2}, v_{3}\}, \ldots, \{v_{n - 1}, v_{n}\}, \{v_{1}, v_{n}\} \}$.  

We may assume that $v_{1} = 1$. 
Then $v_{2} < v_{n - 1}$ since $\{1, v_{2}\} \in E(G)$ but $\{1, v_{n - 1}\}$, $\{v_{2}, v_{n - 1}\} \not\in E(G)$. 
Similarly, we have $v_{n} < v_{3}$ since $\{1, v_{n}\} \in E(G)$ but 
$\{1, v_{3}\}$, $\{v_{3}, v_{n}\} \not\in E(G)$. 
Moreover, we also have $v_{2} < v_{n}$ since $\{v_{n - 1}, v_{n}\} \in E(G)$ but 
$\{v_{2}, v_{n - 1}\}$, $\{v_{2}, v_{n}\} \not\in E(G)$. 
Hence $v_{2} < v_{n} < v_{3}$, but this is a contradiction since $\{v_{2}, v_{3}\} \in E(G)$ but 
$\{v_{2}, v_{n}\}$, $\{v_{3}, v_{n}\} \not\in E(G)$. 
Therefore we have that $C_{n}$ is not weakly closed. 
\end{proof}

\begin{Corollary}
Let $T$ be a tree. 
Then
	\begin{enumerate}[$(1)$]
		\item \cite[Corollary 1.3]{origin} $T$ is closed if and only if $T$ is a path graph. 
		\item $T$ is weakly closed if and only if $T$ is a caterpillar, i.e. a tree for which removing 
		the leaves and incident edges produces a path graph.  
	\end{enumerate}
\end{Corollary}
\begin{proof}
(2) : Let $T$ be a weakly closed tree. 
Then $T$ is bigclaw-free by Proposition \ref{claw}(2). 
Hence $T$ is a caterpillar. 
Conversely, it is easy to see that a caterpillar is weakly closed. 
\end{proof}

Next, we introduce some classes of graphs related to the closed and weakly closed graph. 

\begin{Definition}\label{graphclasses}
	\begin{enumerate}[$(1)$]
		\item Let $P = ([n], <_{P})$ be a partially ordered set. 
		The graph $G(P)$ of $P$ is a graph on the vertex set $[n]$ and $\{i, j\} \in E(G(P))$ with $i < j$ 
		if and only if $i <_{P} j$. 
		A graph $G$ is {\em comparability} if there exists a partially ordered set $P$ such that $G = G(P)$. 
		\item A graph $G$ is {\em complete multipartite} if there exists a partition 
		$V(G) = V_{1} \sqcup \cdots \sqcup V_{r}$ such that $\{u, v\} \not\in E(G)$ if and only if 
		$u, v \in V_{i}$ for some $1 \le i \le r$. 
		\item A graph $G$ is an {\em interval graph} if there exists a set of real intervals 
		$\{I_{v} \mid v \in V(G)\}$ such that $I_{u} \cap I_{v} \neq \emptyset$ if and only if $\{u, v\} \in E(G)$.
		\item A graph $G$ is a {\em proper interval graph} if there exists a set of real intervals 
		$\{I_{v} \mid v \in V(G)\}$ which satisfies that $I_{u} \not\subset I_{v}$ and $I_{v} \not\subset I_{u}$ 
		if $u \neq v$ such that $I_{u} \cap I_{v} \neq \emptyset$ if and only if $\{u, v\} \in E(G)$.
		\item Let $\omega(G)$ be the clique number of $G$ $($i.e. the size of a maximal clique of $G)$ 
		and $\chi(G)$ the chromatic number $($cf. \cite[p.117]{Diestel}$)$ of $G$. 
		A graph $G$ is {\em perfect} if $\omega(H) = \chi(H)$ for all induced subgraph $H$ of $G$. 
	\end{enumerate}
\end{Definition}

\begin{Remark}\label{conditions}
Here we introduce equivalent conditions about graphs appeared in Definition \ref{graphclasses}. 
	\begin{enumerate}[$(1)$]
		\item  Forbidden induced subgraphs of comparability graphs are known 
		$($cf. \cite[p.13]{Mancini}$)$. 
		\item A graph $G$ is complete multipartite if and only if $G$ is $\overline{P}_{3}$-free, 
		where $P_{3}$ is the path graph with $|V(P_{3})| = 3$. 
		\item A graph $G$ is an interval graph if and only if $G$ is chordal and $\overline{G}$ 
		is comparability $($\cite[Theorem 2]{interval}$)$. 
		\item A graph $G$ is a proper interval graph if and only if $G$ is claw-free and interval. 
		\item A graph $G$ is perfect if and only of both $G$ and $\overline{G}$ are 
		$(C_{2n + 3}, n \ge 1)$-free $($strong perfect graph theorem, \cite{CRST}$)$.
	\end{enumerate}
\end{Remark}

Crupi-Rinaldo proved that a graph $G$ is closed if and only if $G$ is a proper interval graph 
(\cite[Theorem 2.4]{CrRi2}). 
Cox-Erskine introduced the notion of narrow graph and they proved that 
a connected graph $G$ is closed if and only if $G$ is chordal, claw-free and narrow 
(\cite[theorem 1.4]{closed1}). 

\begin{Theorem}\label{co-comparability}
Let $G$ be a graph. 
Then the following assertions are equivalent:
\begin{enumerate}[$(1)$]
	\item $G$ is weakly closed.
	\item $G$ is co-comparability, i.e. its complement graph $\overline{G}$ is comparability. 
\end{enumerate}
\end{Theorem}
\begin{proof}
$(1) \Rightarrow (2)$ :  We define the following binary relation $<_{P}$ on $[n]$:
\[
i <_{P} j \iff \{i, j\} \not\in E(G) \ \mbox{and} \ i < j. 
\]
Then $P = ([n], <_{P})$ is a partially ordered set and one has $\overline{G} = G(P)$. 
In fact, we assume that $i <_{P} j$ and $j <_{P} k$. 
Then $\{i, j\}, \{j, k\} \not\in E(G)$ and $i < j < k$. 
Since $G$ is weakly closed, we have $\{i, k\} \not\in E(G)$.  
Hence $i <_{P} k$. 
Thus $<_{P}$ satisfies transitivity. 

$(2) \Rightarrow (1)$ :  Assume that $\overline{G}$ is comparability. 
Then there exists a partially ordered set $P = ([n], <_{P})$ such that 
\[
i <_{P} j \iff \{i, j\} \in E(\overline{G}) \ \mbox{and} \ i < j. 
\]
Assume that $1 \le i < j < k \le n$ and $\{i, k\} \in E(G)$. 
If $\{i, j\}, \{j, k\} \not\in E(G)$, then $i <_{P} j$ and $j <_{P} k$. 
Hence we have $i <_{P} k$, but this is a contradiction since $\{i, k\} \not\in E(\overline{G})$. 
Therefore $G$ is weakly closed. 
\end{proof}

\begin{Corollary}\label{completemultipartite}
Complete multipartite graphs are weakly closed. 
\end{Corollary}
\begin{proof}
Let $G$ be a complete multipartite graph. 
Then $\overline{G}$ is a union of complete graph, hence comparability. 
Thus $G$ is weakly closed by Theorem \ref{co-comparability}. 
\end{proof}

\begin{Corollary}
Interval graphs are weakly closed. 
\end{Corollary}
\begin{proof}
As seen in Remark \ref{conditions} (3), an interval graph
$G$ is chordal and $\overline{G}$ is comparability. 
Hence interval graphs are weakly closed by Theorem \ref{co-comparability}. 
\end{proof}

\begin{Corollary}
Weakly closed graphs are perfect. 
\end{Corollary}
\begin{proof}
Let $G$ be a comparability graph. 
Then $\overline{G}$ is perfect by the Dilworth's theorem \cite{Dilworth}. 
Hence weakly closed graphs are perfect by Theorem \ref{co-comparability}. 
\end{proof}

\section{$F$-purity of binomial edge ideals}

Let $K$ be a perfect field of characteristic $p > 0$ and $S = K[x_{1}, \ldots, x_{n}, y_{1}, \ldots, y_{n}]$ 
the polynomial ring over $K$. 
In this section, we prove that $S/J_{G}$ is $F$-pure if $G$ is weakly closed. 

First, let us set up notations.  
Let $R$ be a reduced Noetherian ring of characteristic $p > 0$ and 
$F:R \to R \ (x \mapsto x^{p})$ the Frobenius map. 
For a non-negative integer $e$, ${}^{e}R$ is the ring $R$ viewed as an 
$R$-module via the $e$-times iterated Frobenius map $F^{e}:R \to R \ (x \mapsto x^{p^{e}})$. 
Now we can identify $F^{e}:R \to {^{e}}R$ with the natural inclusion 
$R \hookrightarrow R^{1/p^{e}}$. 
We say that $R$ is {\em F-finite} if $R^{1/p}$ is a finitely generated $R$-module. 

\begin{Definition}$($\cite{HoR}$)$
Let $R$ be an $F$-finite reduced Noetherian ring of characteristic $p > 0$. 
Then $R$ is said to be {\em F-pure} if the Frobenius map is pure, 
equivalently, the natural inclusion $R \hookrightarrow R^{1/p} \ (x \mapsto (x^p)^{1/p})$ is pure. 
\end{Definition}

It is known that determinantal rings and Stanley-Reisner rings are $F$-pure. 
Since $S/J_{G}$ is reduced (\cite[Corollary 2.2]{origin}), the following question 
is natural: 

\begin{Question}
When is $S/J_{G}$ F-pure ? 
\end{Question}

In \cite{O2}, Ohtani proved that $S/J_{G}$ is $F$-pure if $G$ is complete multipartite. 
Moreover, it is easy to show that $F$-purity of $S/J_{G}$ holds for all closed graphs  
by using Ohtani's technique. 

In this section, we prove the following theorem:

\begin{Theorem}
If $G$ is weakly closed, then $S/J_{G}$ is $F$-pure. 
\end{Theorem} 

Since complete multipartite graphs are weakly closed by Corollary \ref{completemultipartite}, 
this is a generalization of Ohtani's result. 

Before proving this theorem, we prepare some lemmas. 
The first lemma is the graded version of Fedder's criterion. 

\begin{Lemma}$($\cite{Fedder}$)$\label{Fedder}
Let $K$ be a perfect field of characteristic $p > 0$ and $R = K[x_{1}, \ldots, x_{n}]$ the 
polynomial ring over $K$. 
Let $\mathfrak{m} = (x_{1}, \ldots, x_{n})$ be the irrelevant maximal ideal of $R$ and 
let $I \subset \mathfrak{m}$ be a homogeneous ideal of $R$. 
Then $R/I$ is $F$-pure if $I^{[p]} : I \not\subset \mathfrak{m}^{[p]}$, 
where $J^{[p]} = (x^{p} \mid x \in J)$ for an ideal $J \subset S$. 
\end{Lemma}

For an integer sequence $v_{1}, \ldots, v_{n}$ with $\{v_{1}, \ldots, v_{n}\} = \{1, \ldots, n\}$, 
we put
\[
Y_{v_{1}}(v_{1}, v_{2}, \ldots, v_{n})X_{v_{n}} := (Y_{v_{1}} [v_{1}, v_{2}][v_{2}, v_{3}] \cdots [v_{n - 1}, v_{n}] X_{v_{n}})^{p - 1}. 
\]

The second lemma is 

\begin{Lemma}$($\cite[Formula 1, 2]{O2}$)$\label{formula}
If $\{a, b\} \in E(G)$, then
\begin{eqnarray*}
Y_{v_{1}}(v_{1}, \ldots, c, \underline{a}, \underline{b}, d, \ldots, v_{n})X_{v_{n}} &\equiv& Y_{v_{1}}(v_{1}, \ldots, c, \underline{b}, \underline{a}, d, \ldots, v_{n})X_{v_{n}}, \\
Y_{a}(\underline{a}, \underline{b}, c, \ldots, v_{n})X_{v_{n}} &\equiv& Y_{b}(\underline{b}, \underline{a}, c, \ldots, v_{n})X_{v_{n}}, \\
Y_{v_{1}}(v_{1}, \ldots, c, \underline{a}, \underline{b})X_{b} &\equiv& Y_{v_{1}}(v_{1}, \ldots, c, \underline{b}, \underline{a})X_{a}
\end{eqnarray*}
modulo $J^{[p]}_{G}$. 
\end{Lemma}

\begin{proof}[Proof of Theorem 2.3] 
Assume that $G$ is weakly closed. 
By Lemma \ref{Fedder}, it is enough to show that $Y_{1}(1, \ldots, n)X_{n} \in (J^{[p]}_{G} : J_{G}) \setminus \mathfrak{m}^{[p]}$. 

First, we show that $Y_{1}(1, \ldots, n)X_{n} \not\in \mathfrak{m}^{[p]}$. 
Let $>$ be the lexicographic order on $S$ with $x_{1} > \cdots x_{n} > y_{1} > \cdots > y_{n}$. 
Then ${\rm in}_{<} (Y_{1}(1, \ldots, n)X_{n}) = (x_{1} \cdots x_{n}y_{1} \cdots y_{n})^{p - 1}  \not\in \mathfrak{m}^{[p]}$. 
Hence we have $Y_{1}(1, \ldots, n)X_{n} \not\in \mathfrak{m}^{[p]}$. 

Next, we show that $Y_{1}(1, \ldots, n)X_{n} \in J^{[p]}_{G} : J_{G}$. 
Then it is enough to show that $Y_{1}(1, \ldots, n)X_{n} \cdot [i, j] \in J^{[p]}_{G}$ for all $\{i, j\} \in E(G)$. 

Assume that $\{i, k\} \in E(G)$. 
If $k - i = 1$, then it is trivial that $Y_{1}(1, \ldots, n)X_{n} \cdot [i, k] \in J^{[p]}_{G}$. 

Assume that $k - i \ge 2$. 
Now we define 
\begin{eqnarray*}
A&:=&\{j \mid i < j < k, \{j, k\} \in E(G)\}, \\
B&:=&\{j \mid i < j < k, \{j, k\} \not\in E(G)\}. 
\end{eqnarray*}
Since $G$ is weakly closed, we have that for any $s \in A$ and for any $t \in B$, 
$\{s, t\} \in E(G)$ if $s < t$. 
Hence, by using Lemma \ref{formula} repeatedly, we have 
\[
Y_{1}(1, \ldots, n)X_{n} \equiv Y_{1}(1, \ldots, i - 1, i, w_{1}, \ldots, w_{b}, k, v_{1}, \ldots, v_{a}, k +1, \ldots, n)X_{n}
\]
modulo $J^{[p]}_{G}$, where 
\begin{eqnarray*}
A &=& \{v_{1}, \ldots, v_{a}\} \ \mbox{with} \ v_{1} < \cdots < v_{a}, \\
B &=& \{w_{1}, \ldots, w_{b}\} \ \mbox{with} \ w_{1} < \cdots < w_{b}. 
\end{eqnarray*}
Moreover, $\{i, w_{1}\}, \ldots, \{i, w_{b}\} \in E(G)$ since $\{i, k\} \in E(G)$ and 
$\{w_{1}, k\}, \ldots, \{w_{b}, k\} \not\in E(G)$. 
Hence, by using Lemma \ref{formula} again, we have 
\begin{eqnarray*}
& & Y_{1}(1, \ldots, i - 1, i, w_{1}, \ldots, w_{b}, k, v_{1}, \ldots, v_{a}, k +1, \ldots, n)X_{n} \\
&\cong& Y_{1}(1, \ldots, i - 1, w_{1}, \ldots, w_{b}, i, k, v_{1}, \ldots, v_{a}, k +1, \ldots, n)X_{n} \\
&=& (Y_{1}[1, 2] \cdots [i-2, i-1] [i - 1, w_{1}] [w_{1}, w_{2}] \cdots [w_{b - 1}, w_{b}] [w_{b}, i] \\
&\times& \underline{[i, k]} [k, v_{1}] [v_{1}, v_{2}] \cdots [v_{a-1}, v_{a}] [v_{a}, k + 1] [k + 1, k + 2] \cdots [n - 1, n] X_{n})^{p - 1}
\end{eqnarray*}
modulo $J^{[p]}_{G}$. 
Therefore we have that $Y_{1}(1, \ldots, n)X_{n} \in J^{[p]}_{G} : J_{G}$. 
\end{proof}

\begin{Corollary}
We use the same notation us above. 
If $G$ is complete multipartite, then $S/J_{G}$ is $F$-pure. 
\end{Corollary}

In the case of non-weakly closed graph, the $F$-purity of $S/J_{G}$ is quite mysterious. 

\begin{Example}
By using Macaulay2 \cite{Macaulay2}, we can check that $S/J_{C_{5}}$ is not $F$-pure if $p = 2$, but 
is $F$-pure if $p = 3, 5, 7$. 
Moreover, we can also check that $S/J_{C_{6}}$ is not $F$-pure if $p = 2, 3, 5$. 
\end{Example}

\begin{Conjecture}
Assume that the characteristic of the base field is $p = 2$. 
Then $S/J_{G}$ is $F$-pure if and only if $G$ is weakly closed. 
\end{Conjecture}

\begin{Conjecture}
Let $G$ be a graph. 
Then $S/J_{G}$ is $F$-pure for all sufficiently large $p > 0$. 
\end{Conjecture}

\section{Classification of connected weakly closed graphs}

As the end of this paper, 
we classify all connected weakly closed graphs with $5$ and $6$ vertices 
by using Harary's classification of connected graph \cite{Ha}. 

\bigskip

\noindent
{\bf Acknowledgment.}
The author wish to thank Professor Takayuki Hibi for his financial support.  
He also wish to thank Professor Ken-ichi Yoshida and Masahiro Ohtani for their valuable comments. 
He would like to thank Professor Giancarlo Rinaldo for his indication about citation mistake. 
The author was partially supported by JSPS KAKENHI 17K14165. 

The revision of this paper was delayed significantly because of the author's negligence. 
The author is deeply grateful to all researchers who are interested in my work.

\newpage


\begin{xy}
	\ar@{} (0, 0);(50, 0) *\txt{Connected graphs with $5$ vertices (21 items)};
	\ar@{-} (0, 0);(8, 0);
	\ar@{-} (92, 0);(150, 0) = "A";
	\ar@{-} (0, 0);(0, -218) = "G";
	\ar@{-} "A";(150, -218) = "H";
	\ar@{-} "G";"H";
	
	\ar@{} (0, 0);(32, -8) *\cir<2pt>{} = "B1";
	\ar@{-} "B1";(24, -13) *\cir<2pt>{} = "C1";
	\ar@{-} "C1";(27, -20) *\cir<2pt>{} = "D1";
	\ar@{-} "D1";(37, -20) *\cir<2pt>{} = "E1";
	\ar@{-} "E1";(40, -13) *\cir<2pt>{} = "F1";
	\ar@{-} "F1";"B1";
	
	\textcolor{blue}{\ar@{-} (5, -25);(60, -25);}
	\textcolor{blue}{\ar@{} (0, 0);(75, -25) *\txt{Weakly closed};}
	\textcolor{blue}{\ar@{-} (90, -25);(142, -25);}
	\textcolor{blue}{\ar@{-} (1, -25);(1, -215);}
	\textcolor{blue}{\ar@{-} (139, -25);(139, -215);}
	\textcolor{blue}{\ar@{-} (-2, -215);(138, -215);}
	
	\ar@{} (0, 0);(35, -32) *\txt{Complete multipartite};
	\ar@{-} (8, -32);(13, -32);
	\ar@{-} (58, -32);(125, -32) = "AC";
	\ar@{-} (8, -32);(8, -88) = "GC";
	\ar@{-} "AC";(125, -88) = "HC";
	\ar@{-} "GC";"HC";
	
	\ar@{} (0, 0);(23, -38) *\cir<2pt>{} = "B2";
	\ar@{-} "B2";(15, -43) *\cir<2pt>{} = "C2";
	\ar@{-} "C2";(18, -50) *\cir<2pt>{} = "D2";
	\ar@{-} "D2";(28, -50) *\cir<2pt>{} = "E2";
	\ar@{-} "E2";(31, -43) *\cir<2pt>{} = "F2";
	\ar@{-} "F2";"B2";
	\ar@{-} "B2";"D2";
	\ar@{-} "B2";"E2";
	\ar@{-} "C2";"F2";
	
	\ar@{} (0, 0);(54, -38) *\cir<2pt>{} = "B3";
	\ar@{-} "B3";(46, -43) *\cir<2pt>{} = "C3";
	\ar@{-} "C3";(49, -50) *\cir<2pt>{} = "D3";
	\ar@{-} "D3";(59, -50) *\cir<2pt>{} = "E3";
	\ar@{-} "E3";(62, -43) *\cir<2pt>{} = "F3";
	\ar@{-} "B3";"E3";
	\ar@{-} "C3";"F3";

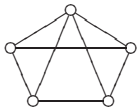
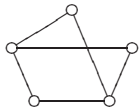
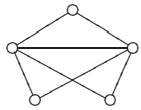
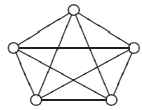
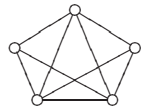
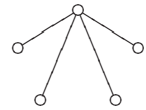
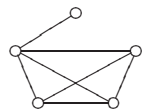
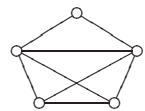
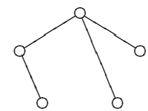
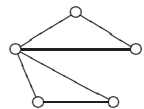
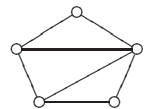
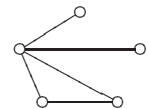
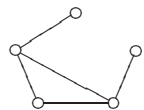
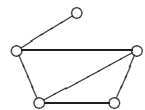
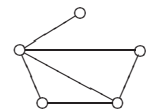
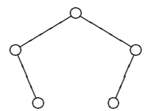
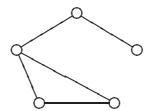
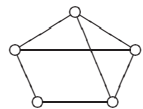
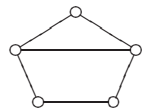
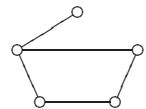
	\textcolor{green}{\ar@{-} (2, -55);(91, -55);}
	\textcolor{green}{\ar@{} (0, 0);(100, -55) *\txt{Chordal};}
	\textcolor{green}{\ar@{-} (109, -55);(129, -55);}
	\textcolor{green}{\ar@{-} (-2, -55);(-2, -190);}
	\textcolor{green}{\ar@{-} (126, -55);(126, -190);}
	\textcolor{green}{\ar@{-} (-5, -190);(125, -190);}

	\ar@{} (0, 0);(15, -68) *\cir<2pt>{} = "B4";
	\ar@{-} "B4";(7, -73) *\cir<2pt>{} = "C4";
	\ar@{-} "C4";(10, -80) *\cir<2pt>{} = "D4";
	\ar@{} "D4";(20, -80) *\cir<2pt>{} = "E4";
	\ar@{-} "E4";(23, -73) *\cir<2pt>{} = "F4";
	\ar@{-} "F4";"B4";
	\ar@{-} "C4";"E4";
	\ar@{-} "C4";"F4";
	\ar@{-} "D4";"F4";
	
	\ar@{} (0, 0);(46, -68) *\cir<2pt>{} = "B5";
	\ar@{-} "B5";(38, -73) *\cir<2pt>{} = "C5";
	\ar@{-} "C5";(41, -80) *\cir<2pt>{} = "D5";
	\ar@{-} "D5";(51, -80) *\cir<2pt>{} = "E5";
	\ar@{-} "E5";(54, -73) *\cir<2pt>{} = "F5";
	\ar@{-} "F5";"B5";
	\ar@{-} "B5";"D5";
	\ar@{-} "B5";"E5";
	\ar@{-} "C5";"E5";
	\ar@{-} "C5";"F5";
	\ar@{-} "D5";"F5";

	\ar@{} (0, 0);(69, -68) *\cir<2pt>{} = "B6";
	\ar@{-} "B6";(61, -73) *\cir<2pt>{} = "C6";
	\ar@{-} "C6";(64, -80) *\cir<2pt>{} = "D6";
	\ar@{-} "D6";(74, -80) *\cir<2pt>{} = "E6";
	\ar@{-} "E6";(77, -73) *\cir<2pt>{} = "F6";
	\ar@{-} "F6";"B6";
	\ar@{-} "B6";"D6";
	\ar@{-} "B6";"E6";
	\ar@{-} "C6";"E6";
	\ar@{} "C6";"F6";
	\ar@{-} "D6";"F6";
	
	\ar@{} (0, 0);(100, -68) *\cir<2pt>{} = "B7";
	\ar@{-} "B7";(92, -73) *\cir<2pt>{} = "C7";
	\ar@{} "C7";(95, -80) *\cir<2pt>{} = "D7";
	\ar@{} "D7";(105, -80) *\cir<2pt>{} = "E7";
	\ar@{} "E7";(108, -73) *\cir<2pt>{} = "F7";
	\ar@{-} "F7";"B7";
	\ar@{-} "B7";"D7";
	\ar@{-} "B7";"E7";
	\ar@{} "C7";"E7";
	\ar@{} "C7";"F7";
	\ar@{} "D7";"F7";
	
	\textcolor{red}{\ar@{-} (31, -62);(48, -62);}
	\textcolor{red}{\ar@{} (0, 0);(56, -62) *\txt{Closed};}
	\textcolor{red}{\ar@{-} (64, -62);(81, -62);}
	\textcolor{red}{\ar@{-} (27, -62);(27, -184);}
	\textcolor{red}{\ar@{-} (78, -62);(78, -184);}
	\textcolor{red}{\ar@{-} (24, -184);(77, -184);}

	\ar@{} (0, 0);(38, -95) *\cir<2pt>{} = "B8";
	\ar@{-} "B8";(30, -100) *\cir<2pt>{} = "C8";
	\ar@{-} "C8";(33, -107) *\cir<2pt>{} = "D8";
	\ar@{-} "D8";(43, -107) *\cir<2pt>{} = "E8";
	\ar@{-} "E8";(46, -100) *\cir<2pt>{} = "F8";
	\ar@{} "F8";"B8";
	\ar@{} "B8";"D8";
	\ar@{} "B8";"E8";
	\ar@{-} "C8";"E8";
	\ar@{-} "C8";"F8";
	\ar@{-} "D8";"F8";

	\ar@{} (0, 0);(61, -95) *\cir<2pt>{} = "B9";
	\ar@{-} "B9";(53, -100) *\cir<2pt>{} = "C9";
	\ar@{-} "C9";(56, -107) *\cir<2pt>{} = "D9";
	\ar@{-} "D9";(66, -107) *\cir<2pt>{} = "E9";
	\ar@{-} "E9";(69, -100) *\cir<2pt>{} = "F9";
	\ar@{-} "F9";"B9";
	\ar@{} "B9";"D9";
	\ar@{} "B9";"E9";
	\ar@{-} "C9";"E9";
	\ar@{-} "C9";"F9";
	\ar@{-} "D9";"F9";
	
	\ar@{} (0, 0);(92, -95) *\cir<2pt>{} = "B10";
	\ar@{-} "B10";(84, -100) *\cir<2pt>{} = "C10";
	\ar@{-} "C10";(87, -107) *\cir<2pt>{} = "D10";
	\ar@{} "D10";(97, -107) *\cir<2pt>{} = "E10";
	\ar@{} "E10";(100, -100) *\cir<2pt>{} = "F10";
	\ar@{-} "F10";"B10";
	\ar@{} "B10";"D10";
	\ar@{-} "B10";"E10";
	\ar@{} "C10";"E10";
	\ar@{} "C10";"F10";
	\ar@{} "D10";"F10";

	\ar@{} (0, 0);(38, -118) *\cir<2pt>{} = "B11";
	\ar@{-} "B11";(30, -123) *\cir<2pt>{} = "C11";
	\ar@{-} "C11";(33, -130) *\cir<2pt>{} = "D11";
	\ar@{-} "D11";(43, -130) *\cir<2pt>{} = "E11";
	\ar@{} "E11";(46, -123) *\cir<2pt>{} = "F11";
	\ar@{-} "F11";"B11";
	\ar@{} "B11";"D11";
	\ar@{} "B11";"E11";
	\ar@{-} "C11";"E11";
	\ar@{-} "C11";"F11";
	\ar@{} "D11";"F11";

	\ar@{} (0, 0);(61, -118) *\cir<2pt>{} = "B12";
	\ar@{-} "B12";(53, -123) *\cir<2pt>{} = "C12";
	\ar@{-} "C12";(56, -130) *\cir<2pt>{} = "D12";
	\ar@{-} "D12";(66, -130) *\cir<2pt>{} = "E12";
	\ar@{-} "E12";(69, -123) *\cir<2pt>{} = "F12";
	\ar@{-} "F12";"B12";
	\ar@{} "B12";"D12";
	\ar@{} "B12";"E12";
	\ar@{} "C12";"E12";
	\ar@{-} "C12";"F12";
	\ar@{-} "D12";"F12";
	
	\ar@{} (0, 0);(92, -118) *\cir<2pt>{} = "B13";
	\ar@{-} "B13";(84, -123) *\cir<2pt>{} = "C13";
	\ar@{-} "C13";(87, -130) *\cir<2pt>{} = "D13";
	\ar@{-} "D13";(97, -130) *\cir<2pt>{} = "E13";
	\ar@{} "E13";(100, -123) *\cir<2pt>{} = "F13";
	\ar@{} "F13";"B13";
	\ar@{} "B13";"D13";
	\ar@{} "B13";"E13";
	\ar@{-} "C13";"E13";
	\ar@{-} "C13";"F13";
	\ar@{} "D13";"F13";

	\ar@{} (0, 0);(38, -141) *\cir<2pt>{} = "B14";
	\ar@{-} "B14";(30, -146) *\cir<2pt>{} = "C14";
	\ar@{-} "C14";(33, -153) *\cir<2pt>{} = "D14";
	\ar@{-} "D14";(43, -153) *\cir<2pt>{} = "E14";
	\ar@{-} "E14";(46, -146) *\cir<2pt>{} = "F14";
	\ar@{} "F14";"B14";
	\ar@{} "B14";"D14";
	\ar@{} "B14";"E14";
	\ar@{-} "C14";"E14";
	\ar@{} "C14";"F14";
	\ar@{} "D14";"F14";

	\ar@{} (0, 0);(61, -141) *\cir<2pt>{} = "B15";
	\ar@{-} "B15";(53, -146) *\cir<2pt>{} = "C15";
	\ar@{-} "C15";(56, -153) *\cir<2pt>{} = "D15";
	\ar@{-} "D15";(66, -153) *\cir<2pt>{} = "E15";
	\ar@{-} "E15";(69, -146) *\cir<2pt>{} = "F15";
	\ar@{} "F15";"B15";
	\ar@{} "B15";"D15";
	\ar@{} "B15";"E15";
	\ar@{} "C15";"E15";
	\ar@{-} "C15";"F15";
	\ar@{-} "D15";"F15";

	\ar@{} (0, 0);(92, -141) *\cir<2pt>{} = "B16";
	\ar@{-} "B16";(84, -146) *\cir<2pt>{} = "C16";
	\ar@{-} "C16";(87, -153) *\cir<2pt>{} = "D16";
	\ar@{-} "D16";(97, -153) *\cir<2pt>{} = "E16";
	\ar@{-} "E16";(100, -146) *\cir<2pt>{} = "F16";
	\ar@{} "F16";"B16";
	\ar@{} "B16";"D16";
	\ar@{} "B16";"E16";
	\ar@{-} "C16";"E16";
	\ar@{-} "C16";"F16";
	\ar@{} "D16";"F16";

	\ar@{} (0, 0);(38, -164) *\cir<2pt>{} = "B17";
	\ar@{-} "B17";(30, -169) *\cir<2pt>{} = "C17";
	\ar@{-} "C17";(33, -176) *\cir<2pt>{} = "D17";
	\ar@{} "D17";(43, -176) *\cir<2pt>{} = "E17";
	\ar@{-} "E17";(46, -169) *\cir<2pt>{} = "F17";
	\ar@{-} "F17";"B17";
	\ar@{} "B17";"D17";
	\ar@{} "B17";"E17";
	\ar@{} "C17";"E17";
	\ar@{} "C17";"F17";
	\ar@{} "D17";"F17";

	\ar@{} (0, 0);(61, -164) *\cir<2pt>{} = "B18";
	\ar@{-} "B18";(53, -169) *\cir<2pt>{} = "C18";
	\ar@{-} "C18";(56, -176) *\cir<2pt>{} = "D18";
	\ar@{-} "D18";(66, -176) *\cir<2pt>{} = "E18";
	\ar@{} "E18";(69, -169) *\cir<2pt>{} = "F18";
	\ar@{-} "F18";"B18";
	\ar@{} "B18";"D18";
	\ar@{} "B18";"E18";
	\ar@{-} "C18";"E18";
	\ar@{} "C18";"F18";
	\ar@{} "D18";"F18";
	
	\ar@{} (0, 0);(10, -198) *\cir<2pt>{} = "B19";
	\ar@{-} "B19";(2, -203) *\cir<2pt>{} = "C19";
	\ar@{-} "C19";(5, -210) *\cir<2pt>{} = "D19";
	\ar@{-} "D19";(15, -210) *\cir<2pt>{} = "E19";
	\ar@{-} "E19";(18, -203) *\cir<2pt>{} = "F19";
	\ar@{-} "F19";"B19";
	\ar@{} "B19";"D19";
	\ar@{-} "B19";"E19";
	\ar@{} "C19";"E19";
	\ar@{-} "C19";"F19";
	\ar@{} "D19";"F19";

	\ar@{} (0, 0);(38, -198) *\cir<2pt>{} = "B20";
	\ar@{-} "B20";(30, -203) *\cir<2pt>{} = "C20";
	\ar@{-} "C20";(33, -210) *\cir<2pt>{} = "D20";
	\ar@{-} "D20";(43, -210) *\cir<2pt>{} = "E20";
	\ar@{-} "E20";(46, -203) *\cir<2pt>{} = "F20";
	\ar@{-} "F20";"B20";
	\ar@{} "B20";"D20";
	\ar@{} "B20";"E20";
	\ar@{} "C20";"E20";
	\ar@{-} "C20";"F20";
	\ar@{} "D20";"F20";

	\ar@{} (0, 0);(66, -198) *\cir<2pt>{} = "B21";
	\ar@{-} "B21";(58, -203) *\cir<2pt>{} = "C21";
	\ar@{-} "C21";(61, -210) *\cir<2pt>{} = "D21";
	\ar@{-} "D21";(71, -210) *\cir<2pt>{} = "E21";
	\ar@{-} "E21";(74, -203) *\cir<2pt>{} = "F21";
	\ar@{} "F21";"B21";
	\ar@{} "B21";"D21";
	\ar@{} "B21";"E21";
	\ar@{} "C21";"E21";
	\ar@{-} "C21";"F21";
	\ar@{} "D21";"F21";
\end{xy}

\newpage

\begin{xy}
	\ar@{} (0, 0);(50, 0) *\txt{Connected graphs with $6$ vertices (112 items)};
	\ar@{-} (0, 0);(8, 0);
	\ar@{-} (92, 0);(150, 0) = "A";
	\ar@{-} (0, 0);(0, -218) = "H";
	\ar@{-} "A";(150, -218) = "I";
	\ar@{-} "H";"I";

	\ar@{} (0, 0);(10, -6) *\cir<2pt>{} = "B1";
	\ar@{-} "B1";(8, -10) *\cir<2pt>{} = "C1";
	\ar@{-} "C1";(10, -14) *\cir<2pt>{} = "D1";
	\ar@{-} "D1";(14, -14) *\cir<2pt>{} = "E1";
	\ar@{-} "E1";(16, -10) *\cir<2pt>{} = "F1";
	\ar@{-} "F1";(14, -6) *\cir<2pt>{} = "G1";
	\ar@{-} "B1";"G1";
	\ar@{} "B1";"D1";
	\ar@{} "B1";"E1";
	\ar@{} "B1";"F1";
	\ar@{} "C1";"E1";
	\ar@{} "C1";"F1";
	\ar@{} "C1";"G1";
	\ar@{} "D1";"F1";
	\ar@{} "D1";"G1";
	\ar@{} "E1";"G1";

	\ar@{} (0, 0);(22, -6) *\cir<2pt>{} = "B2";
	\ar@{-} "B2";(20, -10) *\cir<2pt>{} = "C2";
	\ar@{-} "C2";(22, -14) *\cir<2pt>{} = "D2";
	\ar@{-} "D2";(26, -14) *\cir<2pt>{} = "E2";
	\ar@{-} "E2";(28, -10) *\cir<2pt>{} = "F2";
	\ar@{-} "F2";(26, -6) *\cir<2pt>{} = "G2";
	\ar@{} "B2";"G2";
	\ar@{} "B2";"D2";
	\ar@{} "B2";"E2";
	\ar@{-} "B2";"F2";
	\ar@{} "C2";"E2";
	\ar@{} "C2";"F2";
	\ar@{} "C2";"G2";
	\ar@{} "D2";"F2";
	\ar@{} "D2";"G2";
	\ar@{} "E2";"G2";

	\ar@{} (0, 0);(34, -6) *\cir<2pt>{} = "B3";
	\ar@{-} "B3";(32, -10) *\cir<2pt>{} = "C3";
	\ar@{-} "C3";(34, -14) *\cir<2pt>{} = "D3";
	\ar@{-} "D3";(38, -14) *\cir<2pt>{} = "E3";
	\ar@{-} "E3";(40, -10) *\cir<2pt>{} = "F3";
	\ar@{-} "F3";(38, -6) *\cir<2pt>{} = "G3";
	\ar@{-} "B3";"G3";
	\ar@{} "B3";"D3";
	\ar@{} "B3";"E3";
	\ar@{} "B3";"F3";
	\ar@{} "C3";"E3";
	\ar@{} "C3";"F3";
	\ar@{} "C3";"G3";
	\ar@{} "D3";"F3";
	\ar@{} "D3";"G3";
	\ar@{-} "E3";"G3";

	\ar@{} (0, 0);(46, -6) *\cir<2pt>{} = "B4";
	\ar@{-} "B4";(44, -10) *\cir<2pt>{} = "C4";
	\ar@{-} "C4";(46, -14) *\cir<2pt>{} = "D4";
	\ar@{-} "D4";(50, -14) *\cir<2pt>{} = "E4";
	\ar@{-} "E4";(52, -10) *\cir<2pt>{} = "F4";
	\ar@{-} "F4";(50, -6) *\cir<2pt>{} = "G4";
	\ar@{} "B4";"G4";
	\ar@{} "B4";"D4";
	\ar@{} "B4";"E4";
	\ar@{-} "B4";"F4";
	\ar@{} "C4";"E4";
	\ar@{} "C4";"F4";
	\ar@{-} "C4";"G4";
	\ar@{} "D4";"F4";
	\ar@{} "D4";"G4";
	\ar@{} "E4";"G4";

	\ar@{} (0, 0);(58, -6) *\cir<2pt>{} = "B5";
	\ar@{-} "B5";(56, -10) *\cir<2pt>{} = "C5";
	\ar@{-} "C5";(58, -14) *\cir<2pt>{} = "D5";
	\ar@{-} "D5";(62, -14) *\cir<2pt>{} = "E5";
	\ar@{-} "E5";(64, -10) *\cir<2pt>{} = "F5";
	\ar@{-} "F5";(62, -6) *\cir<2pt>{} = "G5";
	\ar@{-} "B5";"G5";
	\ar@{} "B5";"D5";
	\ar@{} "B5";"E5";
	\ar@{-} "B5";"F5";
	\ar@{} "C5";"E5";
	\ar@{} "C5";"F5";
	\ar@{} "C5";"G5";
	\ar@{} "D5";"F5";
	\ar@{-} "D5";"G5";
	\ar@{} "E5";"G5";

	\ar@{} (0, 0);(70, -6) *\cir<2pt>{} = "B6";
	\ar@{-} "B6";(68, -10) *\cir<2pt>{} = "C6";
	\ar@{-} "C6";(70, -14) *\cir<2pt>{} = "D6";
	\ar@{-} "D6";(74, -14) *\cir<2pt>{} = "E6";
	\ar@{-} "E6";(76, -10) *\cir<2pt>{} = "F6";
	\ar@{-} "F6";(74, -6) *\cir<2pt>{} = "G6";
	\ar@{-} "B6";"G6";
	\ar@{-} "B6";"D6";
	\ar@{} "B6";"E6";
	\ar@{-} "B6";"F6";
	\ar@{} "C6";"E6";
	\ar@{} "C6";"F6";
	\ar@{} "C6";"G6";
	\ar@{} "D6";"F6";
	\ar@{} "D6";"G6";
	\ar@{} "E6";"G6";

	\ar@{} (0, 0);(82, -6) *\cir<2pt>{} = "B7";
	\ar@{-} "B7";(80, -10) *\cir<2pt>{} = "C7";
	\ar@{-} "C7";(82, -14) *\cir<2pt>{} = "D7";
	\ar@{-} "D7";(86, -14) *\cir<2pt>{} = "E7";
	\ar@{-} "E7";(88, -10) *\cir<2pt>{} = "F7";
	\ar@{-} "F7";(86, -6) *\cir<2pt>{} = "G7";
	\ar@{-} "B7";"G7";
	\ar@{} "B7";"D7";
	\ar@{} "B7";"E7";
	\ar@{-} "B7";"F7";
	\ar@{} "C7";"E7";
	\ar@{} "C7";"F7";
	\ar@{-} "C7";"G7";
	\ar@{} "D7";"F7";
	\ar@{} "D7";"G7";
	\ar@{} "E7";"G7";

	\ar@{} (0, 0);(94, -6) *\cir<2pt>{} = "B8";
	\ar@{-} "B8";(92, -10) *\cir<2pt>{} = "C8";
	\ar@{-} "C8";(94, -14) *\cir<2pt>{} = "D8";
	\ar@{-} "D8";(98, -14) *\cir<2pt>{} = "E8";
	\ar@{-} "E8";(100, -10) *\cir<2pt>{} = "F8";
	\ar@{-} "F8";(98, -6) *\cir<2pt>{} = "G8";
	\ar@{-} "B8";"G8";
	\ar@{-} "B8";"D8";
	\ar@{-} "B8";"E8";
	\ar@{} "B8";"F8";
	\ar@{} "C8";"E8";
	\ar@{} "C8";"F8";
	\ar@{-} "C8";"G8";
	\ar@{} "D8";"F8";
	\ar@{} "D8";"G8";
	\ar@{} "E8";"G8";

	\ar@{} (0, 0);(106, -6) *\cir<2pt>{} = "B9";
	\ar@{-} "B9";(104, -10) *\cir<2pt>{} = "C9";
	\ar@{-} "C9";(106, -14) *\cir<2pt>{} = "D9";
	\ar@{-} "D9";(110, -14) *\cir<2pt>{} = "E9";
	\ar@{-} "E9";(112, -10) *\cir<2pt>{} = "F9";
	\ar@{-} "F9";(110, -6) *\cir<2pt>{} = "G9";
	\ar@{-} "B9";"G9";
	\ar@{} "B9";"D9";
	\ar@{-} "B9";"E9";
	\ar@{} "B9";"F9";
	\ar@{-} "C9";"E9";
	\ar@{} "C9";"F9";
	\ar@{} "C9";"G9";
	\ar@{-} "D9";"F9";
	\ar@{} "D9";"G9";
	\ar@{-} "E9";"G9";

	\ar@{} (0, 0);(9, -46) *\cir<2pt>{} = "B111";
	\ar@{-} "B111";(7, -50) *\cir<2pt>{} = "C111";
	\ar@{} "C111";(9, -54) *\cir<2pt>{} = "D111";
	\ar@{-} "D111";(13, -54) *\cir<2pt>{} = "E111";
	\ar@{} "E111";(15, -50) *\cir<2pt>{} = "F111";
	\ar@{-} "F111";(13, -46) *\cir<2pt>{} = "G111";
	\ar@{-} "B111";"G111";
	\ar@{-} "B111";"D111";
	\ar@{} "B111";"E111";
	\ar@{} "B111";"F111";
	\ar@{} "C111";"E111";
	\ar@{} "C111";"F111";
	\ar@{} "C111";"G111";
	\ar@{} "D111";"F111";
	\ar@{-} "D111";"G111";
	\ar@{} "E111";"G111";
	
	\ar@{} (0, 0);(9, -58) *\cir<2pt>{} = "B112";
	\ar@{-} "B112";(7, -62) *\cir<2pt>{} = "C112";
	\ar@{-} "C112";(9, -66) *\cir<2pt>{} = "D112";
	\ar@{-} "D112";(13, -66) *\cir<2pt>{} = "E112";
	\ar@{-} "E112";(15, -62) *\cir<2pt>{} = "F112";
	\ar@{-} "F112";(13, -58) *\cir<2pt>{} = "G112";
	\ar@{-} "B112";"G112";
	\ar@{} "B112";"D112";
	\ar@{} "B112";"E112";
	\ar@{} "B112";"F112";
	\ar@{-} "C112";"E112";
	\ar@{} "C112";"F112";
	\ar@{-} "C112";"G112";
	\ar@{} "D112";"F112";
	\ar@{} "D112";"G112";
	\ar@{-} "E112";"G112";

	\textcolor{blue}{\ar@{-} (19, -20);(97, -20);}
	\textcolor{blue}{\ar@{} (0, 0);(110, -20) *\txt{Weakly closed};}
	\textcolor{blue}{\ar@{-} (123, -20);(142, -20);}
	\textcolor{blue}{\ar@{-} (15, -20);(15, -202);}
	\textcolor{blue}{\ar@{-} (139, -20);(139, -202);}
	\textcolor{blue}{\ar@{-} (12, -202);(138, -202);}
	
	\ar@{} (0, 0);(65, -24) *\txt{Complete multipartite};
	\ar@{-} (38, -24);(42, -24);
	\ar@{-} (88, -24);(105, -24) = "AC";
	\ar@{-} (38, -24);(38, -62) = "GC";
	\ar@{-} "AC";(105, -62) = "HC";
	\ar@{-} "GC";"HC";
	
	\ar@{} (0, 0);(44, -30) *\cir<2pt>{} = "B10";
	\ar@{-} "B10";(42, -34) *\cir<2pt>{} = "C10";
	\ar@{-} "C10";(44, -38) *\cir<2pt>{} = "D10";
	\ar@{-} "D10";(48, -38) *\cir<2pt>{} = "E10";
	\ar@{-} "E10";(50, -34) *\cir<2pt>{} = "F10";
	\ar@{} "F10";(48, -30) *\cir<2pt>{} = "G10";
	\ar@{} "B10";"G10";
	\ar@{} "B10";"D10";
	\ar@{-} "B10";"E10";
	\ar@{} "B10";"F10";
	\ar@{} "C10";"E10";
	\ar@{-} "C10";"F10";
	\ar@{-} "C10";"G10";
	\ar@{} "D10";"F10";
	\ar@{} "D10";"G10";
	\ar@{-} "E10";"G10";
	
	\ar@{} (0, 0);(56, -30) *\cir<2pt>{} = "B11";
	\ar@{-} "B11";(54, -34) *\cir<2pt>{} = "C11";
	\ar@{-} "C11";(56, -38) *\cir<2pt>{} = "D11";
	\ar@{-} "D11";(60, -38) *\cir<2pt>{} = "E11";
	\ar@{-} "E11";(62, -34) *\cir<2pt>{} = "F11";
	\ar@{-} "F11";(60, -30) *\cir<2pt>{} = "G11";
	\ar@{-} "B11";"G11";
	\ar@{} "B11";"D11";
	\ar@{-} "B11";"E11";
	\ar@{} "B11";"F11";
	\ar@{} "C11";"E11";
	\ar@{-} "C11";"F11";
	\ar@{} "C11";"G11";
	\ar@{} "D11";"F11";
	\ar@{-} "D11";"G11";
	\ar@{} "E11";"G11";	
	
	\ar@{} (0, 0);(68, -30) *\cir<2pt>{} = "B12";
	\ar@{-} "B12";(66, -34) *\cir<2pt>{} = "C12";
	\ar@{-} "C12";(68, -38) *\cir<2pt>{} = "D12";
	\ar@{-} "D12";(72, -38) *\cir<2pt>{} = "E12";
	\ar@{-} "E12";(74, -34) *\cir<2pt>{} = "F12";
	\ar@{-} "F12";(72, -30) *\cir<2pt>{} = "G12";
	\ar@{-} "B12";"G12";
	\ar@{} "B12";"D12";
	\ar@{-} "B12";"E12";
	\ar@{} "B12";"F12";
	\ar@{-} "C12";"E12";
	\ar@{-} "C12";"F12";
	\ar@{-} "C12";"G12";
	\ar@{} "D12";"F12";
	\ar@{-} "D12";"G12";
	\ar@{} "E12";"G12";	
	
	\ar@{} (0, 0);(84, -30) *\cir<2pt>{} = "B13";
	\ar@{-} "B13";(82, -34) *\cir<2pt>{} = "C13";
	\ar@{-} "C13";(84, -38) *\cir<2pt>{} = "D13";
	\ar@{-} "D13";(88, -38) *\cir<2pt>{} = "E13";
	\ar@{-} "E13";(90, -34) *\cir<2pt>{} = "F13";
	\ar@{-} "F13";(88, -30) *\cir<2pt>{} = "G13";
	\ar@{-} "B13";"G13";
	\ar@{-} "B13";"D13";
	\ar@{} "B13";"E13";
	\ar@{-} "B13";"F13";
	\ar@{-} "C13";"E13";
	\ar@{} "C13";"F13";
	\ar@{-} "C13";"G13";
	\ar@{-} "D13";"F13";
	\ar@{} "D13";"G13";
	\ar@{-} "E13";"G13";
	
	\ar@{} (0, 0);(96, -30) *\cir<2pt>{} = "B14";
	\ar@{-} "B14";(94, -34) *\cir<2pt>{} = "C14";
	\ar@{-} "C14";(96, -38) *\cir<2pt>{} = "D14";
	\ar@{-} "D14";(100, -38) *\cir<2pt>{} = "E14";
	\ar@{-} "E14";(102, -34) *\cir<2pt>{} = "F14";
	\ar@{-} "F14";(100, -30) *\cir<2pt>{} = "G14";
	\ar@{-} "B14";"G14";
	\ar@{} "B14";"D14";
	\ar@{-} "B14";"E14";
	\ar@{-} "B14";"F14";
	\ar@{-} "C14";"E14";
	\ar@{-} "C14";"F14";
	\ar@{-} "C14";"G14";
	\ar@{-} "D14";"F14";
	\ar@{-} "D14";"G14";
	\ar@{} "E14";"G14";	
	
	\ar@{} (0, 0);(44, -50) *\cir<2pt>{} = "B15";
	\ar@{-} "B15";(42, -54) *\cir<2pt>{} = "C15";
	\ar@{-} "C15";(44, -58) *\cir<2pt>{} = "D15";
	\ar@{-} "D15";(48, -58) *\cir<2pt>{} = "E15";
	\ar@{-} "E15";(50, -54) *\cir<2pt>{} = "F15";
	\ar@{-} "F15";(48, -50) *\cir<2pt>{} = "G15";
	\ar@{-} "B15";"G15";
	\ar@{-} "B15";"D15";
	\ar@{-} "B15";"E15";
	\ar@{-} "B15";"F15";
	\ar@{} "C15";"E15";
	\ar@{-} "C15";"F15";
	\ar@{} "C15";"G15";
	\ar@{-} "D15";"F15";
	\ar@{-} "D15";"G15";
	\ar@{} "E15";"G15";		
	
	\ar@{} (0, 0);(56, -50) *\cir<2pt>{} = "B16";
	\ar@{-} "B16";(54, -54) *\cir<2pt>{} = "C16";
	\ar@{-} "C16";(56, -58) *\cir<2pt>{} = "D16";
	\ar@{} "D16";(60, -58) *\cir<2pt>{} = "E16";
	\ar@{-} "E16";(62, -54) *\cir<2pt>{} = "F16";
	\ar@{-} "F16";(60, -50) *\cir<2pt>{} = "G16";
	\ar@{} "B16";"G16";
	\ar@{} "B16";"D16";
	\ar@{} "B16";"E16";
	\ar@{-} "B16";"F16";
	\ar@{-} "C16";"E16";
	\ar@{-} "C16";"F16";
	\ar@{-} "C16";"G16";
	\ar@{-} "D16";"F16";
	\ar@{} "D16";"G16";
	\ar@{} "E16";"G16";	
	
	\ar@{} (0, 0);(68, -50) *\cir<2pt>{} = "B17";
	\ar@{} "B17";(66, -54) *\cir<2pt>{} = "C17";
	\ar@{} "C17";(68, -58) *\cir<2pt>{} = "D17";
	\ar@{} "D17";(72, -58) *\cir<2pt>{} = "E17";
	\ar@{} "E17";(74, -54) *\cir<2pt>{} = "F17";
	\ar@{-} "F17";(72, -50) *\cir<2pt>{} = "G17";
	\ar@{-} "B17";"G17";
	\ar@{} "B17";"D17";
	\ar@{} "B17";"E17";
	\ar@{} "B17";"F17";
	\ar@{} "C17";"E17";
	\ar@{} "C17";"F17";
	\ar@{-} "C17";"G17";
	\ar@{} "D17";"F17";
	\ar@{-} "D17";"G17";
	\ar@{-} "E17";"G17";	
	
	\ar@{} (0, 0);(84, -50) *\cir<2pt>{} = "B18";
	\ar@{-} "B18";(82, -54) *\cir<2pt>{} = "C18";
	\ar@{-} "C18";(84, -58) *\cir<2pt>{} = "D18";
	\ar@{-} "D18";(88, -58) *\cir<2pt>{} = "E18";
	\ar@{-} "E18";(90, -54) *\cir<2pt>{} = "F18";
	\ar@{-} "F18";(88, -50) *\cir<2pt>{} = "G18";
	\ar@{-} "B18";"G18";
	\ar@{-} "B18";"D18";
	\ar@{-} "B18";"E18";
	\ar@{-} "B18";"F18";
	\ar@{-} "C18";"E18";
	\ar@{-} "C18";"F18";
	\ar@{-} "C18";"G18";
	\ar@{-} "D18";"F18";
	\ar@{-} "D18";"G18";
	\ar@{-} "E18";"G18";	
	
	\ar@{} (0, 0);(96, -50) *\cir<2pt>{} = "B19";
	\ar@{-} "B19";(94, -54) *\cir<2pt>{} = "C19";
	\ar@{-} "C19";(96, -58) *\cir<2pt>{} = "D19";
	\ar@{-} "D19";(100, -58) *\cir<2pt>{} = "E19";
	\ar@{-} "E19";(102, -54) *\cir<2pt>{} = "F19";
	\ar@{-} "F19";(100, -50) *\cir<2pt>{} = "G19";
	\ar@{-} "B19";"G19";
	\ar@{} "B19";"D19";
	\ar@{-} "B19";"E19";
	\ar@{-} "B19";"F19";
	\ar@{-} "C19";"E19";
	\ar@{-} "C19";"F19";
	\ar@{-} "C19";"G19";
	\ar@{-} "D19";"F19";
	\ar@{-} "D19";"G19";
	\ar@{-} "E19";"G19";	
	
	\ar@{} (0, 0);(84, -66) *\cir<2pt>{} = "B20";
	\ar@{-} "B20";(82, -70) *\cir<2pt>{} = "C20";
	\ar@{-} "C20";(84, -74) *\cir<2pt>{} = "D20";
	\ar@{-} "D20";(88, -74) *\cir<2pt>{} = "E20";
	\ar@{-} "E20";(90, -70) *\cir<2pt>{} = "F20";
	\ar@{} "F20";(88, -66) *\cir<2pt>{} = "G20";
	\ar@{-} "B20";"G20";
	\ar@{} "B20";"D20";
	\ar@{} "B20";"E20";
	\ar@{} "B20";"F20";
	\ar@{} "C20";"E20";
	\ar@{} "C20";"F20";
	\ar@{} "C20";"G20";
	\ar@{} "D20";"F20";
	\ar@{} "D20";"G20";
	\ar@{} "E20";"G20";	
	
	\ar@{} (0, 0);(96, -66) *\cir<2pt>{} = "B21";
	\ar@{-} "B21";(94, -70) *\cir<2pt>{} = "C21";
	\ar@{-} "C21";(96, -74) *\cir<2pt>{} = "D21";
	\ar@{-} "D21";(100, -74) *\cir<2pt>{} = "E21";
	\ar@{-} "E21";(102, -70) *\cir<2pt>{} = "F21";
	\ar@{-} "F21";(100, -66) *\cir<2pt>{} = "G21";
	\ar@{} "B21";"G21";
	\ar@{-} "B21";"D21";
	\ar@{-} "B21";"E21";
	\ar@{} "B21";"F21";
	\ar@{} "C21";"E21";
	\ar@{} "C21";"F21";
	\ar@{} "C21";"G21";
	\ar@{} "D21";"F21";
	\ar@{} "D21";"G21";
	\ar@{-} "E21";"G21";	
	
	\ar@{} (0, 0);(108, -66) *\cir<2pt>{} = "B22";
	\ar@{-} "B22";(106, -70) *\cir<2pt>{} = "C22";
	\ar@{-} "C22";(108, -74) *\cir<2pt>{} = "D22";
	\ar@{-} "D22";(112, -74) *\cir<2pt>{} = "E22";
	\ar@{-} "E22";(114, -70) *\cir<2pt>{} = "F22";
	\ar@{-} "F22";(112, -66) *\cir<2pt>{} = "G22";
	\ar@{-} "B22";"G22";
	\ar@{-} "B22";"D22";
	\ar@{-} "B22";"E22";
	\ar@{} "B22";"F22";
	\ar@{} "C22";"E22";
	\ar@{} "C22";"F22";
	\ar@{} "C22";"G22";
	\ar@{} "D22";"F22";
	\ar@{-} "D22";"G22";
	\ar@{-} "E22";"G22";	
	
	\ar@{} (0, 0);(120, -66) *\cir<2pt>{} = "B23";
	\ar@{-} "B23";(118, -70) *\cir<2pt>{} = "C23";
	\ar@{-} "C23";(120, -74) *\cir<2pt>{} = "D23";
	\ar@{-} "D23";(124, -74) *\cir<2pt>{} = "E23";
	\ar@{-} "E23";(126, -70) *\cir<2pt>{} = "F23";
	\ar@{-} "F23";(124, -66) *\cir<2pt>{} = "G23";
	\ar@{} "B23";"G23";
	\ar@{-} "B23";"D23";
	\ar@{} "B23";"E23";
	\ar@{} "B23";"F23";
	\ar@{} "C23";"E23";
	\ar@{} "C23";"F23";
	\ar@{} "C23";"G23";
	\ar@{} "D23";"F23";
	\ar@{} "D23";"G23";
	\ar@{} "E23";"G23";	
	
	\ar@{} (0, 0);(84, -78) *\cir<2pt>{} = "B24";
	\ar@{-} "B24";(82, -82) *\cir<2pt>{} = "C24";
	\ar@{-} "C24";(84, -86) *\cir<2pt>{} = "D24";
	\ar@{-} "D24";(88, -86) *\cir<2pt>{} = "E24";
	\ar@{} "E24";(90, -82) *\cir<2pt>{} = "F24";
	\ar@{-} "F24";(88, -78) *\cir<2pt>{} = "G24";
	\ar@{-} "B24";"G24";
	\ar@{-} "B24";"D24";
	\ar@{-} "B24";"E24";
	\ar@{} "B24";"F24";
	\ar@{} "C24";"E24";
	\ar@{} "C24";"F24";
	\ar@{} "C24";"G24";
	\ar@{} "D24";"F24";
	\ar@{} "D24";"G24";
	\ar@{-} "E24";"G24";	

	\ar@{} (0, 0);(96, -78) *\cir<2pt>{} = "B25";
	\ar@{-} "B25";(94, -82) *\cir<2pt>{} = "C25";
	\ar@{-} "C25";(96, -86) *\cir<2pt>{} = "D25";
	\ar@{-} "D25";(100, -86) *\cir<2pt>{} = "E25";
	\ar@{-} "E25";(102, -82) *\cir<2pt>{} = "F25";
	\ar@{-} "F25";(100, -78) *\cir<2pt>{} = "G25";
	\ar@{-} "B25";"G25";
	\ar@{-} "B25";"D25";
	\ar@{-} "B25";"E25";
	\ar@{} "B25";"F25";
	\ar@{-} "C25";"E25";
	\ar@{} "C25";"F25";
	\ar@{} "C25";"G25";
	\ar@{} "D25";"F25";
	\ar@{} "D25";"G25";
	\ar@{-} "E25";"G25";	
	
	\ar@{} (0, 0);(108, -78) *\cir<2pt>{} = "B26";
	\ar@{-} "B26";(106, -82) *\cir<2pt>{} = "C26";
	\ar@{-} "C26";(108, -86) *\cir<2pt>{} = "D26";
	\ar@{-} "D26";(112, -86) *\cir<2pt>{} = "E26";
	\ar@{} "E26";(114, -82) *\cir<2pt>{} = "F26";
	\ar@{-} "F26";(112, -78) *\cir<2pt>{} = "G26";
	\ar@{-} "B26";"G26";
	\ar@{-} "B26";"D26";
	\ar@{} "B26";"E26";
	\ar@{} "B26";"F26";
	\ar@{} "C26";"E26";
	\ar@{} "C26";"F26";
	\ar@{} "C26";"G26";
	\ar@{} "D26";"F26";
	\ar@{} "D26";"G26";
	\ar@{} "E26";"G26";	
	
	\ar@{} (0, 0);(120, -78) *\cir<2pt>{} = "B27";
	\ar@{} "B27";(118, -82) *\cir<2pt>{} = "C27";
	\ar@{-} "C27";(120, -86) *\cir<2pt>{} = "D27";
	\ar@{-} "D27";(124, -86) *\cir<2pt>{} = "E27";
	\ar@{-} "E27";(126, -82) *\cir<2pt>{} = "F27";
	\ar@{} "F27";(124, -78) *\cir<2pt>{} = "G27";
	\ar@{-} "B27";"G27";
	\ar@{-} "B27";"D27";
	\ar@{-} "B27";"E27";
	\ar@{} "B27";"F27";
	\ar@{} "C27";"E27";
	\ar@{} "C27";"F27";
	\ar@{} "C27";"G27";
	\ar@{} "D27";"F27";
	\ar@{-} "D27";"G27";
	\ar@{-} "E27";"G27";	
	
	\ar@{} (0, 0);(84, -90) *\cir<2pt>{} = "B28";
	\ar@{-} "B28";(82, -94) *\cir<2pt>{} = "C28";
	\ar@{-} "C28";(84, -98) *\cir<2pt>{} = "D28";
	\ar@{-} "D28";(88, -98) *\cir<2pt>{} = "E28";
	\ar@{-} "E28";(90, -94) *\cir<2pt>{} = "F28";
	\ar@{} "F28";(88, -90) *\cir<2pt>{} = "G28";
	\ar@{-} "B28";"G28";
	\ar@{-} "B28";"D28";
	\ar@{-} "B28";"E28";
	\ar@{} "B28";"F28";
	\ar@{-} "C28";"E28";
	\ar@{} "C28";"F28";
	\ar@{-} "C28";"G28";
	\ar@{} "D28";"F28";
	\ar@{-} "D28";"G28";
	\ar@{-} "E28";"G28";	
	
	\ar@{} (0, 0);(96, -90) *\cir<2pt>{} = "B29";
	\ar@{-} "B29";(94, -94) *\cir<2pt>{} = "C29";
	\ar@{-} "C29";(96, -98) *\cir<2pt>{} = "D29";
	\ar@{-} "D29";(100, -98) *\cir<2pt>{} = "E29";
	\ar@{} "E29";(102, -94) *\cir<2pt>{} = "F29";
	\ar@{-} "F29";(100, -90) *\cir<2pt>{} = "G29";
	\ar@{} "B29";"G29";
	\ar@{-} "B29";"D29";
	\ar@{-} "B29";"E29";
	\ar@{} "B29";"F29";
	\ar@{} "C29";"E29";
	\ar@{} "C29";"F29";
	\ar@{} "C29";"G29";
	\ar@{} "D29";"F29";
	\ar@{} "D29";"G29";
	\ar@{-} "E29";"G29";	
	
	\ar@{} (0, 0);(108, -90) *\cir<2pt>{} = "B30";
	\ar@{-} "B30";(106, -94) *\cir<2pt>{} = "C30";
	\ar@{-} "C30";(108, -98) *\cir<2pt>{} = "D30";
	\ar@{-} "D30";(112, -98) *\cir<2pt>{} = "E30";
	\ar@{-} "E30";(114, -94) *\cir<2pt>{} = "F30";
	\ar@{} "F30";(112, -90) *\cir<2pt>{} = "G30";
	\ar@{-} "B30";"G30";
	\ar@{-} "B30";"D30";
	\ar@{-} "B30";"E30";
	\ar@{} "B30";"F30";
	\ar@{} "C30";"E30";
	\ar@{} "C30";"F30";
	\ar@{} "C30";"G30";
	\ar@{} "D30";"F30";
	\ar@{-} "D30";"G30";
	\ar@{-} "E30";"G30";	
	
	\ar@{} (0, 0);(120, -90) *\cir<2pt>{} = "B31";
	\ar@{-} "B31";(118, -94) *\cir<2pt>{} = "C31";
	\ar@{-} "C31";(120, -98) *\cir<2pt>{} = "D31";
	\ar@{-} "D31";(124, -98) *\cir<2pt>{} = "E31";
	\ar@{-} "E31";(126, -94) *\cir<2pt>{} = "F31";
	\ar@{-} "F31";(124, -90) *\cir<2pt>{} = "G31";
	\ar@{-} "B31";"G31";
	\ar@{} "B31";"D31";
	\ar@{} "B31";"E31";
	\ar@{-} "B31";"F31";
	\ar@{-} "C31";"E31";
	\ar@{-} "C31";"F31";
	\ar@{-} "C31";"G31";
	\ar@{-} "D31";"F31";
	\ar@{} "D31";"G31";
	\ar@{} "E31";"G31";	

	\ar@{} (0, 0);(84, -102) *\cir<2pt>{} = "B32";
	\ar@{-} "B32";(82, -106) *\cir<2pt>{} = "C32";
	\ar@{-} "C32";(84, -110) *\cir<2pt>{} = "D32";
	\ar@{-} "D32";(88, -110) *\cir<2pt>{} = "E32";
	\ar@{-} "E32";(90, -106) *\cir<2pt>{} = "F32";
	\ar@{} "F32";(88, -102) *\cir<2pt>{} = "G32";
	\ar@{} "B32";"G32";
	\ar@{-} "B32";"D32";
	\ar@{} "B32";"E32";
	\ar@{} "B32";"F32";
	\ar@{} "C32";"E32";
	\ar@{} "C32";"F32";
	\ar@{} "C32";"G32";
	\ar@{} "D32";"F32";
	\ar@{-} "D32";"G32";
	\ar@{-} "E32";"G32";	

	\ar@{} (0, 0);(96, -102) *\cir<2pt>{} = "B33";
	\ar@{-} "B33";(94, -106) *\cir<2pt>{} = "C33";
	\ar@{-} "C33";(96, -110) *\cir<2pt>{} = "D33";
	\ar@{-} "D33";(100, -110) *\cir<2pt>{} = "E33";
	\ar@{-} "E33";(102, -106) *\cir<2pt>{} = "F33";
	\ar@{-} "F33";(100, -102) *\cir<2pt>{} = "G33";
	\ar@{-} "B33";"G33";
	\ar@{-} "B33";"D33";
	\ar@{} "B33";"E33";
	\ar@{} "B33";"F33";
	\ar@{} "C33";"E33";
	\ar@{} "C33";"F33";
	\ar@{} "C33";"G33";
	\ar@{} "D33";"F33";
	\ar@{-} "D33";"G33";
	\ar@{-} "E33";"G33";	

	\ar@{} (0, 0);(108, -102) *\cir<2pt>{} = "B34";
	\ar@{-} "B34";(106, -106) *\cir<2pt>{} = "C34";
	\ar@{-} "C34";(108, -110) *\cir<2pt>{} = "D34";
	\ar@{-} "D34";(112, -110) *\cir<2pt>{} = "E34";
	\ar@{-} "E34";(114, -106) *\cir<2pt>{} = "F34";
	\ar@{-} "F34";(112, -102) *\cir<2pt>{} = "G34";
	\ar@{} "B34";"G34";
	\ar@{-} "B34";"D34";
	\ar@{-} "B34";"E34";
	\ar@{-} "B34";"F34";
	\ar@{} "C34";"E34";
	\ar@{-} "C34";"F34";
	\ar@{-} "C34";"G34";
	\ar@{-} "D34";"F34";
	\ar@{} "D34";"G34";
	\ar@{} "E34";"G34";	

	\ar@{} (0, 0);(120, -102) *\cir<2pt>{} = "B35";
	\ar@{-} "B35";(118, -106) *\cir<2pt>{} = "C35";
	\ar@{} "C35";(120, -110) *\cir<2pt>{} = "D35";
	\ar@{-} "D35";(124, -110) *\cir<2pt>{} = "E35";
	\ar@{-} "E35";(126, -106) *\cir<2pt>{} = "F35";
	\ar@{} "F35";(124, -102) *\cir<2pt>{} = "G35";
	\ar@{-} "B35";"G35";
	\ar@{-} "B35";"D35";
	\ar@{} "B35";"E35";
	\ar@{} "B35";"F35";
	\ar@{} "C35";"E35";
	\ar@{} "C35";"F35";
	\ar@{} "C35";"G35";
	\ar@{} "D35";"F35";
	\ar@{-} "D35";"G35";
	\ar@{-} "E35";"G35";	

	\ar@{} (0, 0);(84, -114) *\cir<2pt>{} = "B36";
	\ar@{-} "B36";(82, -118) *\cir<2pt>{} = "C36";
	\ar@{-} "C36";(84, -122) *\cir<2pt>{} = "D36";
	\ar@{-} "D36";(88, -122) *\cir<2pt>{} = "E36";
	\ar@{-} "E36";(90, -118) *\cir<2pt>{} = "F36";
	\ar@{-} "F36";(88, -114) *\cir<2pt>{} = "G36";
	\ar@{} "B36";"G36";
	\ar@{-} "B36";"D36";
	\ar@{-} "B36";"E36";
	\ar@{} "B36";"F36";
	\ar@{-} "C36";"E36";
	\ar@{} "C36";"F36";
	\ar@{} "C36";"G36";
	\ar@{} "D36";"F36";
	\ar@{} "D36";"G36";
	\ar@{-} "E36";"G36";	

	\ar@{} (0, 0);(96, -114) *\cir<2pt>{} = "B37";
	\ar@{-} "B37";(94, -118) *\cir<2pt>{} = "C37";
	\ar@{-} "C37";(96, -122) *\cir<2pt>{} = "D37";
	\ar@{-} "D37";(100, -122) *\cir<2pt>{} = "E37";
	\ar@{-} "E37";(102, -118) *\cir<2pt>{} = "F37";
	\ar@{-} "F37";(100, -114) *\cir<2pt>{} = "G37";
	\ar@{-} "B37";"G37";
	\ar@{} "B37";"D37";
	\ar@{-} "B37";"E37";
	\ar@{-} "B37";"F37";
	\ar@{-} "C37";"E37";
	\ar@{-} "C37";"F37";
	\ar@{-} "C37";"G37";
	\ar@{-} "D37";"F37";
	\ar@{} "D37";"G37";
	\ar@{} "E37";"G37";	

	\ar@{} (0, 0);(108, -114) *\cir<2pt>{} = "B38";
	\ar@{-} "B38";(106, -118) *\cir<2pt>{} = "C38";
	\ar@{-} "C38";(108, -122) *\cir<2pt>{} = "D38";
	\ar@{-} "D38";(112, -122) *\cir<2pt>{} = "E38";
	\ar@{-} "E38";(114, -118) *\cir<2pt>{} = "F38";
	\ar@{-} "F38";(112, -114) *\cir<2pt>{} = "G38";
	\ar@{} "B38";"G38";
	\ar@{-} "B38";"D38";
	\ar@{} "B38";"E38";
	\ar@{} "B38";"F38";
	\ar@{} "C38";"E38";
	\ar@{} "C38";"F38";
	\ar@{} "C38";"G38";
	\ar@{} "D38";"F38";
	\ar@{} "D38";"G38";
	\ar@{-} "E38";"G38";	

	\ar@{} (0, 0);(120, -114) *\cir<2pt>{} = "B39";
	\ar@{-} "B39";(118, -118) *\cir<2pt>{} = "C39";
	\ar@{-} "C39";(120, -122) *\cir<2pt>{} = "D39";
	\ar@{-} "D39";(124, -122) *\cir<2pt>{} = "E39";
	\ar@{-} "E39";(126, -118) *\cir<2pt>{} = "F39";
	\ar@{} "F39";(124, -114) *\cir<2pt>{} = "G39";
	\ar@{-} "B39";"G39";
	\ar@{-} "B39";"D39";
	\ar@{} "B39";"E39";
	\ar@{} "B39";"F39";
	\ar@{} "C39";"E39";
	\ar@{} "C39";"F39";
	\ar@{-} "C39";"G39";
	\ar@{} "D39";"F39";
	\ar@{-} "D39";"G39";
	\ar@{-} "E39";"G39";	

	\ar@{} (0, 0);(84, -126) *\cir<2pt>{} = "B40";
	\ar@{-} "B40";(82, -130) *\cir<2pt>{} = "C40";
	\ar@{-} "C40";(84, -134) *\cir<2pt>{} = "D40";
	\ar@{-} "D40";(88, -134) *\cir<2pt>{} = "E40";
	\ar@{-} "E40";(90, -130) *\cir<2pt>{} = "F40";
	\ar@{-} "F40";(88, -126) *\cir<2pt>{} = "G40";
	\ar@{-} "B40";"G40";
	\ar@{-} "B40";"D40";
	\ar@{-} "B40";"E40";
	\ar@{-} "B40";"F40";
	\ar@{-} "C40";"E40";
	\ar@{-} "C40";"F40";	
	\ar@{} "C40";"G40";
	\ar@{-} "D40";"F40";
	\ar@{} "D40";"G40";
	\ar@{} "E40";"G40";	

	\ar@{} (0, 0);(96, -126) *\cir<2pt>{} = "B41";
	\ar@{-} "B41";(94, -130) *\cir<2pt>{} = "C41";
	\ar@{-} "C41";(96, -134) *\cir<2pt>{} = "D41";
	\ar@{-} "D41";(100, -134) *\cir<2pt>{} = "E41";
	\ar@{} "E41";(102, -130) *\cir<2pt>{} = "F41";
	\ar@{-} "F41";(100, -126) *\cir<2pt>{} = "G41";
	\ar@{-} "B41";"G41";
	\ar@{-} "B41";"D41";
	\ar@{-} "B41";"E41";
	\ar@{} "B41";"F41";
	\ar@{-} "C41";"E41";
	\ar@{} "C41";"F41";	
	\ar@{} "C41";"G41";
	\ar@{} "D41";"F41";
	\ar@{} "D41";"G41";
	\ar@{} "E41";"G41";	

	\ar@{} (0, 0);(108, -126) *\cir<2pt>{} = "B42";
	\ar@{-} "B42";(106, -130) *\cir<2pt>{} = "C42";
	\ar@{-} "C42";(108, -134) *\cir<2pt>{} = "D42";
	\ar@{-} "D42";(112, -134) *\cir<2pt>{} = "E42";
	\ar@{-} "E42";(114, -130) *\cir<2pt>{} = "F42";
	\ar@{} "F42";(112, -126) *\cir<2pt>{} = "G42";
	\ar@{-} "B42";"G42";
	\ar@{-} "B42";"D42";
	\ar@{} "B42";"E42";
	\ar@{} "B42";"F42";
	\ar@{-} "C42";"E42";
	\ar@{} "C42";"F42";	
	\ar@{-} "C42";"G42";
	\ar@{} "D42";"F42";
	\ar@{-} "D42";"G42";
	\ar@{-} "E42";"G42";	

	\ar@{} (0, 0);(120, -126) *\cir<2pt>{} = "B43";
	\ar@{-} "B43";(118, -130) *\cir<2pt>{} = "C43";
	\ar@{-} "C43";(120, -134) *\cir<2pt>{} = "D43";
	\ar@{-} "D43";(124, -134) *\cir<2pt>{} = "E43";
	\ar@{-} "E43";(126, -130) *\cir<2pt>{} = "F43";
	\ar@{-} "F43";(124, -126) *\cir<2pt>{} = "G43";
	\ar@{-} "B43";"G43";
	\ar@{} "B43";"D43";
	\ar@{-} "B43";"E43";
	\ar@{} "B43";"F43";
	\ar@{-} "C43";"E43";
	\ar@{-} "C43";"F43";	
	\ar@{-} "C43";"G43";
	\ar@{-} "D43";"F43";
	\ar@{-} "D43";"G43";
	\ar@{-} "E43";"G43";	

	\ar@{} (0, 0);(18, -66) *\cir<2pt>{} = "B44";
	\ar@{-} "B44";(16, -70) *\cir<2pt>{} = "C44";
	\ar@{-} "C44";(18, -74) *\cir<2pt>{} = "D44";
	\ar@{-} "D44";(22, -74) *\cir<2pt>{} = "E44";
	\ar@{} "E44";(24, -70) *\cir<2pt>{} = "F44";
	\ar@{} "F44";(22, -66) *\cir<2pt>{} = "G44";
	\ar@{} "B44";"G44";
	\ar@{} "B44";"D44";
	\ar@{} "B44";"E44";
	\ar@{} "B44";"F44";
	\ar@{} "C44";"E44";
	\ar@{} "C44";"F44";
	\ar@{} "C44";"G44";
	\ar@{-} "D44";"F44";
	\ar@{-} "D44";"G44";
	\ar@{} "E44";"G44";	

	\ar@{} (0, 0);(30, -66) *\cir<2pt>{} = "B45";
	\ar@{-} "B45";(28, -70) *\cir<2pt>{} = "C45";
	\ar@{} "C45";(30, -74) *\cir<2pt>{} = "D45";
	\ar@{-} "D45";(34, -74) *\cir<2pt>{} = "E45";
	\ar@{} "E45";(36, -70) *\cir<2pt>{} = "F45";
	\ar@{-} "F45";(34, -66) *\cir<2pt>{} = "G45";
	\ar@{-} "B45";"G45";
	\ar@{-} "B45";"D45";
	\ar@{} "B45";"E45";
	\ar@{} "B45";"F45";
	\ar@{} "C45";"E45";
	\ar@{} "C45";"F45";
	\ar@{} "C45";"G45";
	\ar@{} "D45";"F45";
	\ar@{} "D45";"G45";
	\ar@{} "E45";"G45";	
	
	\ar@{} (0, 0);(42, -66) *\cir<2pt>{} = "B46";
	\ar@{-} "B46";(40, -70) *\cir<2pt>{} = "C46";
	\ar@{} "C46";(42, -74) *\cir<2pt>{} = "D46";
	\ar@{} "D46";(46, -74) *\cir<2pt>{} = "E46";
	\ar@{} "E46";(48, -70) *\cir<2pt>{} = "F46";
	\ar@{-} "F46";(46, -66) *\cir<2pt>{} = "G46";
	\ar@{-} "B46";"G46";
	\ar@{-} "B46";"D46";
	\ar@{} "B46";"E46";
	\ar@{} "B46";"F46";
	\ar@{} "C46";"E46";
	\ar@{} "C46";"F46";
	\ar@{} "C46";"G46";
	\ar@{} "D46";"F46";
	\ar@{} "D46";"G46";
	\ar@{-} "E46";"G46";	
	
	\ar@{} (0, 0);(54, -66) *\cir<2pt>{} = "B47";
	\ar@{-} "B47";(52, -70) *\cir<2pt>{} = "C47";
	\ar@{-} "C47";(54, -74) *\cir<2pt>{} = "D47";
	\ar@{-} "D47";(58, -74) *\cir<2pt>{} = "E47";
	\ar@{-} "E47";(60, -70) *\cir<2pt>{} = "F47";
	\ar@{} "F47";(58, -66) *\cir<2pt>{} = "G47";
	\ar@{} "B47";"G47";
	\ar@{} "B47";"D47";
	\ar@{} "B47";"E47";
	\ar@{} "B47";"F47";
	\ar@{} "C47";"E47";
	\ar@{} "C47";"F47";
	\ar@{} "C47";"G47";
	\ar@{} "D47";"F47";
	\ar@{} "D47";"G47";
	\ar@{-} "E47";"G47";	
	
	\ar@{} (0, 0);(66, -66) *\cir<2pt>{} = "B48";
	\ar@{-} "B48";(64, -70) *\cir<2pt>{} = "C48";
	\ar@{-} "C48";(66, -74) *\cir<2pt>{} = "D48";
	\ar@{} "D48";(70, -74) *\cir<2pt>{} = "E48";
	\ar@{} "E48";(72, -70) *\cir<2pt>{} = "F48";
	\ar@{} "F48";(70, -66) *\cir<2pt>{} = "G48";
	\ar@{-} "B48";"G48";
	\ar@{-} "B48";"D48";
	\ar@{-} "B48";"E48";
	\ar@{-} "B48";"F48";
	\ar@{} "C48";"E48";
	\ar@{} "C48";"F48";
	\ar@{} "C48";"G48";
	\ar@{} "D48";"F48";
	\ar@{} "D48";"G48";
	\ar@{} "E48";"G48";	
	
	\ar@{} (0, 0);(18, -78) *\cir<2pt>{} = "B49";
	\ar@{-} "B49";(16, -82) *\cir<2pt>{} = "C49";
	\ar@{-} "C49";(18, -86) *\cir<2pt>{} = "D49";
	\ar@{} "D49";(22, -86) *\cir<2pt>{} = "E49";
	\ar@{} "E49";(24, -82) *\cir<2pt>{} = "F49";
	\ar@{-} "F49";(22, -78) *\cir<2pt>{} = "G49";
	\ar@{-} "B49";"G49";
	\ar@{-} "B49";"D49";
	\ar@{-} "B49";"E49";
	\ar@{} "B49";"F49";
	\ar@{} "C49";"E49";
	\ar@{} "C49";"F49";
	\ar@{} "C49";"G49";
	\ar@{} "D49";"F49";
	\ar@{} "D49";"G49";
	\ar@{} "E49";"G49";	

	\ar@{} (0, 0);(30, -78) *\cir<2pt>{} = "B50";
	\ar@{-} "B50";(28, -82) *\cir<2pt>{} = "C50";
	\ar@{-} "C50";(30, -86) *\cir<2pt>{} = "D50";
	\ar@{} "D50";(34, -86) *\cir<2pt>{} = "E50";
	\ar@{} "E50";(36, -82) *\cir<2pt>{} = "F50";
	\ar@{-} "F50";(34, -78) *\cir<2pt>{} = "G50";
	\ar@{-} "B50";"G50";
	\ar@{-} "B50";"D50";
	\ar@{} "B50";"E50";
	\ar@{} "B50";"F50";
	\ar@{} "C50";"E50";
	\ar@{} "C50";"F50";
	\ar@{} "C50";"G50";
	\ar@{} "D50";"F50";
	\ar@{} "D50";"G50";
	\ar@{-} "E50";"G50";	

	\ar@{} (0, 0);(42, -78) *\cir<2pt>{} = "B51";
	\ar@{-} "B51";(40, -82) *\cir<2pt>{} = "C51";
	\ar@{-} "C51";(42, -86) *\cir<2pt>{} = "D51";
	\ar@{-} "D51";(46, -86) *\cir<2pt>{} = "E51";
	\ar@{} "E51";(48, -82) *\cir<2pt>{} = "F51";
	\ar@{} "F51";(46, -78) *\cir<2pt>{} = "G51";
	\ar@{-} "B51";"G51";
	\ar@{-} "B51";"D51";
	\ar@{} "B51";"E51";
	\ar@{-} "B51";"F51";
	\ar@{} "C51";"E51";
	\ar@{} "C51";"F51";
	\ar@{} "C51";"G51";
	\ar@{} "D51";"F51";
	\ar@{} "D51";"G51";
	\ar@{} "E51";"G51";	

	\ar@{} (0, 0);(54, -78) *\cir<2pt>{} = "B52";
	\ar@{-} "B52";(52, -82) *\cir<2pt>{} = "C52";
	\ar@{} "C52";(54, -86) *\cir<2pt>{} = "D52";
	\ar@{-} "D52";(58, -86) *\cir<2pt>{} = "E52";
	\ar@{} "E52";(60, -82) *\cir<2pt>{} = "F52";
	\ar@{-} "F52";(58, -78) *\cir<2pt>{} = "G52";
	\ar@{-} "B52";"G52";
	\ar@{-} "B52";"D52";
	\ar@{} "B52";"E52";
	\ar@{} "B52";"F52";
	\ar@{} "C52";"E52";
	\ar@{} "C52";"F52";
	\ar@{} "C52";"G52";
	\ar@{} "D52";"F52";
	\ar@{-} "D52";"G52";
	\ar@{-} "E52";"G52";	

	\ar@{} (0, 0);(66, -78) *\cir<2pt>{} = "B53";
	\ar@{-} "B53";(64, -82) *\cir<2pt>{} = "C53";
	\ar@{-} "C53";(66, -86) *\cir<2pt>{} = "D53";
	\ar@{-} "D53";(70, -86) *\cir<2pt>{} = "E53";
	\ar@{} "E53";(72, -82) *\cir<2pt>{} = "F53";
	\ar@{-} "F53";(70, -78) *\cir<2pt>{} = "G53";
	\ar@{-} "B53";"G53";
	\ar@{-} "B53";"D53";
	\ar@{-} "B53";"E53";
	\ar@{} "B53";"F53";
	\ar@{} "C53";"E53";
	\ar@{} "C53";"F53";
	\ar@{} "C53";"G53";
	\ar@{} "D53";"F53";
	\ar@{} "D53";"G53";
	\ar@{} "E53";"G53";	

	\ar@{} (0, 0);(18, -90) *\cir<2pt>{} = "B54";
	\ar@{-} "B54";(16, -94) *\cir<2pt>{} = "C54";
	\ar@{-} "C54";(18, -98) *\cir<2pt>{} = "D54";
	\ar@{-} "D54";(22, -98) *\cir<2pt>{} = "E54";
	\ar@{-} "E54";(24, -94) *\cir<2pt>{} = "F54";
	\ar@{} "F54";(22, -90) *\cir<2pt>{} = "G54";
	\ar@{} "B54";"G54";
	\ar@{-} "B54";"D54";
	\ar@{} "B54";"E54";
	\ar@{} "B54";"F54";
	\ar@{} "C54";"E54";
	\ar@{} "C54";"F54";
	\ar@{} "C54";"G54";
	\ar@{-} "D54";"F54";
	\ar@{-} "D54";"G54";
	\ar@{} "E54";"G54";	

	\ar@{} (0, 0);(30, -90) *\cir<2pt>{} = "B55";
	\ar@{-} "B55";(28, -94) *\cir<2pt>{} = "C55";
	\ar@{-} "C55";(30, -98) *\cir<2pt>{} = "D55";
	\ar@{-} "D55";(34, -98) *\cir<2pt>{} = "E55";
	\ar@{-} "E55";(36, -94) *\cir<2pt>{} = "F55";
	\ar@{} "F55";(34, -90) *\cir<2pt>{} = "G55";
	\ar@{} "B55";"G55";
	\ar@{} "B55";"D55";
	\ar@{-} "B55";"E55";
	\ar@{} "B55";"F55";
	\ar@{-} "C55";"E55";
	\ar@{} "C55";"F55";
	\ar@{} "C55";"G55";
	\ar@{} "D55";"F55";
	\ar@{} "D55";"G55";
	\ar@{-} "E55";"G55";	

	\ar@{} (0, 0);(42, -90) *\cir<2pt>{} = "B56";
	\ar@{} "B56";(40, -94) *\cir<2pt>{} = "C56";
	\ar@{-} "C56";(42, -98) *\cir<2pt>{} = "D56";
	\ar@{-} "D56";(46, -98) *\cir<2pt>{} = "E56";
	\ar@{} "E56";(48, -94) *\cir<2pt>{} = "F56";
	\ar@{-} "F56";(46, -90) *\cir<2pt>{} = "G56";
	\ar@{-} "B56";"G56";
	\ar@{-} "B56";"D56";
	\ar@{} "B56";"E56";
	\ar@{} "B56";"F56";
	\ar@{} "C56";"E56";
	\ar@{} "C56";"F56";
	\ar@{} "C56";"G56";
	\ar@{} "D56";"F56";
	\ar@{-} "D56";"G56";
	\ar@{-} "E56";"G56";	

	\ar@{} (0, 0);(54, -90) *\cir<2pt>{} = "B57";
	\ar@{-} "B57";(52, -94) *\cir<2pt>{} = "C57";
	\ar@{-} "C57";(54, -98) *\cir<2pt>{} = "D57";
	\ar@{} "D57";(58, -98) *\cir<2pt>{} = "E57";
	\ar@{} "E57";(60, -94) *\cir<2pt>{} = "F57";
	\ar@{-} "F57";(58, -90) *\cir<2pt>{} = "G57";
	\ar@{-} "B57";"G57";
	\ar@{-} "B57";"D57";
	\ar@{} "B57";"E57";
	\ar@{} "B57";"F57";
	\ar@{} "C57";"E57";
	\ar@{} "C57";"F57";
	\ar@{} "C57";"G57";
	\ar@{} "D57";"F57";
	\ar@{-} "D57";"G57";
	\ar@{-} "E57";"G57";	

	\ar@{} (0, 0);(66, -90) *\cir<2pt>{} = "B58";
	\ar@{-} "B58";(64, -94) *\cir<2pt>{} = "C58";
	\ar@{-} "C58";(66, -98) *\cir<2pt>{} = "D58";
	\ar@{-} "D58";(70, -98) *\cir<2pt>{} = "E58";
	\ar@{-} "E58";(72, -94) *\cir<2pt>{} = "F58";
	\ar@{-} "F58";(70, -90) *\cir<2pt>{} = "G58";
	\ar@{} "B58";"G58";
	\ar@{} "B58";"D58";
	\ar@{-} "B58";"E58";
	\ar@{} "B58";"F58";
	\ar@{-} "C58";"E58";
	\ar@{} "C58";"F58";
	\ar@{} "C58";"G58";
	\ar@{} "D58";"F58";
	\ar@{} "D58";"G58";
	\ar@{-} "E58";"G58";	

	\ar@{} (0, 0);(18, -102) *\cir<2pt>{} = "B59";
	\ar@{-} "B59";(16, -106) *\cir<2pt>{} = "C59";
	\ar@{-} "C59";(18, -110) *\cir<2pt>{} = "D59";
	\ar@{-} "D59";(22, -110) *\cir<2pt>{} = "E59";
	\ar@{-} "E59";(24, -106) *\cir<2pt>{} = "F59";
	\ar@{} "F59";(22, -102) *\cir<2pt>{} = "G59";
	\ar@{-} "B59";"G59";
	\ar@{-} "B59";"D59";
	\ar@{-} "B59";"E59";
	\ar@{} "B59";"F59";
	\ar@{} "C59";"E59";
	\ar@{} "C59";"F59";
	\ar@{} "C59";"G59";
	\ar@{} "D59";"F59";
	\ar@{} "D59";"G59";
	\ar@{-} "E59";"G59";	

	\ar@{} (0, 0);(30, -102) *\cir<2pt>{} = "B60";
	\ar@{-} "B60";(28, -106) *\cir<2pt>{} = "C60";
	\ar@{-} "C60";(30, -110) *\cir<2pt>{} = "D60";
	\ar@{-} "D60";(34, -110) *\cir<2pt>{} = "E60";
	\ar@{-} "E60";(36, -106) *\cir<2pt>{} = "F60";
	\ar@{} "F60";(34, -102) *\cir<2pt>{} = "G60";
	\ar@{-} "B60";"G60";
	\ar@{} "B60";"D60";
	\ar@{-} "B60";"E60";
	\ar@{} "B60";"F60";
	\ar@{-} "C60";"E60";
	\ar@{} "C60";"F60";
	\ar@{} "C60";"G60";
	\ar@{} "D60";"F60";
	\ar@{} "D60";"G60";
	\ar@{-} "E60";"G60";	
	
	\ar@{} (0, 0);(42, -102) *\cir<2pt>{} = "B61";
	\ar@{-} "B61";(40, -106) *\cir<2pt>{} = "C61";
	\ar@{-} "C61";(42, -110) *\cir<2pt>{} = "D61";
	\ar@{-} "D61";(46, -110) *\cir<2pt>{} = "E61";
	\ar@{-} "E61";(48, -106) *\cir<2pt>{} = "F61";
	\ar@{} "F61";(46, -102) *\cir<2pt>{} = "G61";
	\ar@{} "B61";"G61";
	\ar@{} "B61";"D61";
	\ar@{-} "B61";"E61";
	\ar@{} "B61";"F61";
	\ar@{-} "C61";"E61";
	\ar@{} "C61";"F61";
	\ar@{-} "C61";"G61";
	\ar@{} "D61";"F61";
	\ar@{} "D61";"G61";
	\ar@{-} "E61";"G61";	

	\ar@{} (0, 0);(54, -102) *\cir<2pt>{} = "B62";
	\ar@{-} "B62";(52, -106) *\cir<2pt>{} = "C62";
	\ar@{-} "C62";(54, -110) *\cir<2pt>{} = "D62";
	\ar@{-} "D62";(58, -110) *\cir<2pt>{} = "E62";
	\ar@{-} "E62";(60, -106) *\cir<2pt>{} = "F62";
	\ar@{} "F62";(58, -102) *\cir<2pt>{} = "G62";
	\ar@{} "B62";"G62";
	\ar@{-} "B62";"D62";
	\ar@{-} "B62";"E62";
	\ar@{} "B62";"F62";
	\ar@{-} "C62";"E62";
	\ar@{} "C62";"F62";
	\ar@{} "C62";"G62";
	\ar@{} "D62";"F62";
	\ar@{} "D62";"G62";
	\ar@{-} "E62";"G62";	

	\ar@{} (0, 0);(66, -102) *\cir<2pt>{} = "B63";
	\ar@{-} "B63";(64, -106) *\cir<2pt>{} = "C63";
	\ar@{-} "C63";(66, -110) *\cir<2pt>{} = "D63";
	\ar@{-} "D63";(70, -110) *\cir<2pt>{} = "E63";
	\ar@{-} "E63";(72, -106) *\cir<2pt>{} = "F63";
	\ar@{} "F63";(70, -102) *\cir<2pt>{} = "G63";
	\ar@{-} "B63";"G63";
	\ar@{-} "B63";"D63";
	\ar@{-} "B63";"E63";
	\ar@{} "B63";"F63";
	\ar@{} "C63";"E63";
	\ar@{} "C63";"F63";
	\ar@{} "C63";"G63";
	\ar@{} "D63";"F63";
	\ar@{-} "D63";"G63";
	\ar@{} "E63";"G63";	

	\ar@{} (0, 0);(18, -114) *\cir<2pt>{} = "B64";
	\ar@{-} "B64";(16, -118) *\cir<2pt>{} = "C64";
	\ar@{-} "C64";(18, -122) *\cir<2pt>{} = "D64";
	\ar@{-} "D64";(22, -122) *\cir<2pt>{} = "E64";
	\ar@{-} "E64";(24, -118) *\cir<2pt>{} = "F64";
	\ar@{-} "F64";(22, -114) *\cir<2pt>{} = "G64";
	\ar@{-} "B64";"G64";
	\ar@{-} "B64";"D64";
	\ar@{-} "B64";"E64";
	\ar@{-} "B64";"F64";
	\ar@{} "C64";"E64";
	\ar@{} "C64";"F64";
	\ar@{} "C64";"G64";
	\ar@{} "D64";"F64";
	\ar@{} "D64";"G64";
	\ar@{} "E64";"G64";	

	\ar@{} (0, 0);(30, -114) *\cir<2pt>{} = "B65";
	\ar@{-} "B65";(28, -118) *\cir<2pt>{} = "C65";
	\ar@{-} "C65";(30, -122) *\cir<2pt>{} = "D65";
	\ar@{-} "D65";(34, -122) *\cir<2pt>{} = "E65";
	\ar@{-} "E65";(36, -118) *\cir<2pt>{} = "F65";
	\ar@{} "F65";(34, -114) *\cir<2pt>{} = "G65";
	\ar@{-} "B65";"G65";
	\ar@{-} "B65";"D65";
	\ar@{-} "B65";"E65";
	\ar@{} "B65";"F65";
	\ar@{-} "C65";"E65";
	\ar@{} "C65";"F65";
	\ar@{} "C65";"G65";
	\ar@{} "D65";"F65";
	\ar@{} "D65";"G65";
	\ar@{-} "E65";"G65";	
	
	\ar@{} (0, 0);(42, -114) *\cir<2pt>{} = "B66";
	\ar@{-} "B66";(40, -118) *\cir<2pt>{} = "C66";
	\ar@{-} "C66";(42, -122) *\cir<2pt>{} = "D66";
	\ar@{-} "D66";(46, -122) *\cir<2pt>{} = "E66";
	\ar@{-} "E66";(48, -118) *\cir<2pt>{} = "F66";
	\ar@{-} "F66";(46, -114) *\cir<2pt>{} = "G66";
	\ar@{} "B66";"G66";
	\ar@{} "B66";"D66";
	\ar@{} "B66";"E66";
	\ar@{-} "B66";"F66";
	\ar@{} "C66";"E66";
	\ar@{-} "C66";"F66";
	\ar@{-} "C66";"G66";
	\ar@{-} "D66";"F66";
	\ar@{} "D66";"G66";
	\ar@{} "E66";"G66";	
	
	\ar@{} (0, 0);(54, -114) *\cir<2pt>{} = "B67";
	\ar@{-} "B67";(52, -118) *\cir<2pt>{} = "C67";
	\ar@{-} "C67";(54, -122) *\cir<2pt>{} = "D67";
	\ar@{-} "D67";(58, -122) *\cir<2pt>{} = "E67";
	\ar@{-} "E67";(60, -118) *\cir<2pt>{} = "F67";
	\ar@{-} "F67";(58, -114) *\cir<2pt>{} = "G67";
	\ar@{-} "B67";"G67";
	\ar@{} "B67";"D67";
	\ar@{-} "B67";"E67";
	\ar@{} "B67";"F67";
	\ar@{-} "C67";"E67";
	\ar@{} "C67";"F67";
	\ar@{-} "C67";"G67";
	\ar@{} "D67";"F67";
	\ar@{} "D67";"G67";
	\ar@{-} "E67";"G67";	

	\ar@{} (0, 0);(66, -114) *\cir<2pt>{} = "B68";
	\ar@{-} "B68";(64, -118) *\cir<2pt>{} = "C68";
	\ar@{-} "C68";(66, -122) *\cir<2pt>{} = "D68";
	\ar@{-} "D68";(70, -122) *\cir<2pt>{} = "E68";
	\ar@{-} "E68";(72, -118) *\cir<2pt>{} = "F68";
	\ar@{-} "F68";(70, -114) *\cir<2pt>{} = "G68";
	\ar@{} "B68";"G68";
	\ar@{} "B68";"D68";
	\ar@{} "B68";"E68";
	\ar@{-} "B68";"F68";
	\ar@{-} "C68";"E68";
	\ar@{-} "C68";"F68";
	\ar@{-} "C68";"G68";
	\ar@{-} "D68";"F68";
	\ar@{} "D68";"G68";
	\ar@{} "E68";"G68";	

	\ar@{} (0, 0);(18, -126) *\cir<2pt>{} = "B69";
	\ar@{-} "B69";(16, -130) *\cir<2pt>{} = "C69";
	\ar@{-} "C69";(18, -134) *\cir<2pt>{} = "D69";
	\ar@{-} "D69";(22, -134) *\cir<2pt>{} = "E69";
	\ar@{-} "E69";(24, -130) *\cir<2pt>{} = "F69";
	\ar@{} "F69";(22, -126) *\cir<2pt>{} = "G69";
	\ar@{-} "B69";"G69";
	\ar@{-} "B69";"D69";
	\ar@{-} "B69";"E69";
	\ar@{} "B69";"F69";
	\ar@{-} "C69";"E69";
	\ar@{} "C69";"F69";
	\ar@{} "C69";"G69";
	\ar@{} "D69";"F69";
	\ar@{-} "D69";"G69";
	\ar@{-} "E69";"G69";	
	
	\ar@{} (0, 0);(30, -126) *\cir<2pt>{} = "B70";
	\ar@{-} "B70";(28, -130) *\cir<2pt>{} = "C70";
	\ar@{-} "C70";(30, -134) *\cir<2pt>{} = "D70";
	\ar@{-} "D70";(34, -134) *\cir<2pt>{} = "E70";
	\ar@{-} "E70";(36, -130) *\cir<2pt>{} = "F70";
	\ar@{-} "F70";(34, -126) *\cir<2pt>{} = "G70";
	\ar@{} "B70";"G70";
	\ar@{} "B70";"D70";
	\ar@{-} "B70";"E70";
	\ar@{-} "B70";"F70";
	\ar@{-} "C70";"E70";
	\ar@{-} "C70";"F70";
	\ar@{-} "C70";"G70";
	\ar@{-} "D70";"F70";
	\ar@{} "D70";"G70";
	\ar@{} "E70";"G70";	

	\ar@{} (0, 0);(18, -150) *\cir<2pt>{} = "B71";
	\ar@{-} "B71";(16, -154) *\cir<2pt>{} = "C71";
	\ar@{-} "C71";(18, -158) *\cir<2pt>{} = "D71";
	\ar@{-} "D71";(22, -158) *\cir<2pt>{} = "E71";
	\ar@{-} "E71";(24, -154) *\cir<2pt>{} = "F71";
	\ar@{} "F71";(22, -150) *\cir<2pt>{} = "G71";
	\ar@{-} "B71";"G71";
	\ar@{} "B71";"D71";
	\ar@{-} "B71";"E71";
	\ar@{} "B71";"F71";
	\ar@{} "C71";"E71";
	\ar@{} "C71";"F71";
	\ar@{} "C71";"G71";
	\ar@{} "D71";"F71";
	\ar@{} "D71";"G71";
	\ar@{} "E71";"G71";	

	\ar@{} (0, 0);(30, -150) *\cir<2pt>{} = "B72";
	\ar@{-} "B72";(28, -154) *\cir<2pt>{} = "C72";
	\ar@{-} "C72";(30, -158) *\cir<2pt>{} = "D72";
	\ar@{-} "D72";(34, -158) *\cir<2pt>{} = "E72";
	\ar@{-} "E72";(36, -154) *\cir<2pt>{} = "F72";
	\ar@{-} "F72";(34, -150) *\cir<2pt>{} = "G72";
	\ar@{} "B72";"G72";
	\ar@{} "B72";"D72";
	\ar@{-} "B72";"E72";
	\ar@{} "B72";"F72";
	\ar@{} "C72";"E72";
	\ar@{} "C72";"F72";
	\ar@{} "C72";"G72";
	\ar@{} "D72";"F72";
	\ar@{} "D72";"G72";
	\ar@{} "E72";"G72";	
	
	\ar@{} (0, 0);(42, -150) *\cir<2pt>{} = "B73";
	\ar@{-} "B73";(40, -154) *\cir<2pt>{} = "C73";
	\ar@{-} "C73";(42, -158) *\cir<2pt>{} = "D73";
	\ar@{-} "D73";(46, -158) *\cir<2pt>{} = "E73";
	\ar@{} "E73";(48, -154) *\cir<2pt>{} = "F73";
	\ar@{} "F73";(46, -150) *\cir<2pt>{} = "G73";
	\ar@{-} "B73";"G73";
	\ar@{} "B73";"D73";
	\ar@{-} "B73";"E73";
	\ar@{-} "B73";"F73";
	\ar@{} "C73";"E73";
	\ar@{} "C73";"F73";
	\ar@{} "C73";"G73";
	\ar@{} "D73";"F73";
	\ar@{} "D73";"G73";
	\ar@{} "E73";"G73";	
	
	\ar@{} (0, 0);(54, -150) *\cir<2pt>{} = "B74";
	\ar@{} "B74";(52, -154) *\cir<2pt>{} = "C74";
	\ar@{-} "C74";(54, -158) *\cir<2pt>{} = "D74";
	\ar@{-} "D74";(58, -158) *\cir<2pt>{} = "E74";
	\ar@{} "E74";(60, -154) *\cir<2pt>{} = "F74";
	\ar@{-} "F74";(58, -150) *\cir<2pt>{} = "G74";
	\ar@{-} "B74";"G74";
	\ar@{-} "B74";"D74";
	\ar@{} "B74";"E74";
	\ar@{} "B74";"F74";
	\ar@{} "C74";"E74";
	\ar@{} "C74";"F74";
	\ar@{} "C74";"G74";
	\ar@{} "D74";"F74";
	\ar@{} "D74";"G74";
	\ar@{-} "E74";"G74";	

	\ar@{} (0, 0);(66, -150) *\cir<2pt>{} = "B75";
	\ar@{-} "B75";(64, -154) *\cir<2pt>{} = "C75";
	\ar@{-} "C75";(66, -158) *\cir<2pt>{} = "D75";
	\ar@{-} "D75";(70, -158) *\cir<2pt>{} = "E75";
	\ar@{-} "E75";(72, -154) *\cir<2pt>{} = "F75";
	\ar@{-} "F75";(70, -150) *\cir<2pt>{} = "G75";
	\ar@{-} "B75";"G75";
	\ar@{} "B75";"D75";
	\ar@{} "B75";"E75";
	\ar@{} "B75";"F75";
	\ar@{} "C75";"E75";
	\ar@{} "C75";"F75";
	\ar@{} "C75";"G75";
	\ar@{} "D75";"F75";
	\ar@{-} "D75";"G75";
	\ar@{} "E75";"G75";	

	\ar@{} (0, 0);(78, -150) *\cir<2pt>{} = "B76";
	\ar@{-} "B76";(76, -154) *\cir<2pt>{} = "C76";
	\ar@{-} "C76";(78, -158) *\cir<2pt>{} = "D76";
	\ar@{-} "D76";(82, -158) *\cir<2pt>{} = "E76";
	\ar@{-} "E76";(84, -154) *\cir<2pt>{} = "F76";
	\ar@{-} "F76";(82, -150) *\cir<2pt>{} = "G76";
	\ar@{} "B76";"G76";
	\ar@{} "B76";"D76";
	\ar@{-} "B76";"E76";
	\ar@{} "B76";"F76";
	\ar@{} "C76";"E76";
	\ar@{} "C76";"F76";
	\ar@{} "C76";"G76";
	\ar@{} "D76";"F76";
	\ar@{} "D76";"G76";
	\ar@{-} "E76";"G76";	

	\ar@{} (0, 0);(90, -150) *\cir<2pt>{} = "B77";
	\ar@{} "B77";(88, -154) *\cir<2pt>{} = "C77";
	\ar@{-} "C77";(90, -158) *\cir<2pt>{} = "D77";
	\ar@{-} "D77";(94, -158) *\cir<2pt>{} = "E77";
	\ar@{-} "E77";(96, -154) *\cir<2pt>{} = "F77";
	\ar@{} "F77";(94, -150) *\cir<2pt>{} = "G77";
	\ar@{-} "B77";"G77";
	\ar@{-} "B77";"D77";
	\ar@{} "B77";"E77";
	\ar@{-} "B77";"F77";
	\ar@{} "C77";"E77";
	\ar@{} "C77";"F77";
	\ar@{} "C77";"G77";
	\ar@{} "D77";"F77";
	\ar@{} "D77";"G77";
	\ar@{-} "E77";"G77";	

	\ar@{} (0, 0);(102, -150) *\cir<2pt>{} = "B78";
	\ar@{-} "B78";(100, -154) *\cir<2pt>{} = "C78";
	\ar@{-} "C78";(102, -158) *\cir<2pt>{} = "D78";
	\ar@{-} "D78";(106, -158) *\cir<2pt>{} = "E78";
	\ar@{-} "E78";(108, -154) *\cir<2pt>{} = "F78";
	\ar@{} "F78";(106, -150) *\cir<2pt>{} = "G78";
	\ar@{-} "B78";"G78";
	\ar@{-} "B78";"D78";
	\ar@{} "B78";"E78";
	\ar@{} "B78";"F78";
	\ar@{} "C78";"E78";
	\ar@{} "C78";"F78";
	\ar@{} "C78";"G78";
	\ar@{} "D78";"F78";
	\ar@{} "D78";"G78";
	\ar@{-} "E78";"G78";	

	\ar@{} (0, 0);(114, -150) *\cir<2pt>{} = "B79";
	\ar@{-} "B79";(112, -154) *\cir<2pt>{} = "C79";
	\ar@{-} "C79";(114, -158) *\cir<2pt>{} = "D79";
	\ar@{-} "D79";(118, -158) *\cir<2pt>{} = "E79";
	\ar@{} "E79";(120, -154) *\cir<2pt>{} = "F79";
	\ar@{-} "F79";(118, -150) *\cir<2pt>{} = "G79";
	\ar@{-} "B79";"G79";
	\ar@{} "B79";"D79";
	\ar@{-} "B79";"E79";
	\ar@{} "B79";"F79";
	\ar@{} "C79";"E79";
	\ar@{} "C79";"F79";
	\ar@{} "C79";"G79";
	\ar@{} "D79";"F79";
	\ar@{} "D79";"G79";
	\ar@{-} "E79";"G79";	

	\ar@{} (0, 0);(126, -150) *\cir<2pt>{} = "B80";
	\ar@{-} "B80";(124, -154) *\cir<2pt>{} = "C80";
	\ar@{} "C80";(126, -158) *\cir<2pt>{} = "D80";
	\ar@{-} "D80";(130, -158) *\cir<2pt>{} = "E80";
	\ar@{-} "E80";(132, -154) *\cir<2pt>{} = "F80";
	\ar@{} "F80";(130, -150) *\cir<2pt>{} = "G80";
	\ar@{-} "B80";"G80";
	\ar@{-} "B80";"D80";
	\ar@{} "B80";"E80";
	\ar@{-} "B80";"F80";
	\ar@{} "C80";"E80";
	\ar@{} "C80";"F80";
	\ar@{} "C80";"G80";
	\ar@{} "D80";"F80";
	\ar@{} "D80";"G80";
	\ar@{-} "E80";"G80";	

	\ar@{} (0, 0);(18, -162) *\cir<2pt>{} = "B81";
	\ar@{-} "B81";(16, -166) *\cir<2pt>{} = "C81";
	\ar@{-} "C81";(18, -170) *\cir<2pt>{} = "D81";
	\ar@{-} "D81";(22, -170) *\cir<2pt>{} = "E81";
	\ar@{-} "E81";(24, -166) *\cir<2pt>{} = "F81";
	\ar@{} "F81";(22, -162) *\cir<2pt>{} = "G81";
	\ar@{-} "B81";"G81";
	\ar@{} "B81";"D81";
	\ar@{-} "B81";"E81";
	\ar@{} "B81";"F81";
	\ar@{} "C81";"E81";
	\ar@{} "C81";"F81";
	\ar@{} "C81";"G81";
	\ar@{} "D81";"F81";
	\ar@{} "D81";"G81";
	\ar@{-} "E81";"G81";	

	\ar@{} (0, 0);(30, -162) *\cir<2pt>{} = "B82";
	\ar@{-} "B82";(28, -166) *\cir<2pt>{} = "C82";
	\ar@{-} "C82";(30, -170) *\cir<2pt>{} = "D82";
	\ar@{-} "D82";(34, -170) *\cir<2pt>{} = "E82";
	\ar@{-} "E82";(36, -166) *\cir<2pt>{} = "F82";
	\ar@{-} "F82";(34, -162) *\cir<2pt>{} = "G82";
	\ar@{-} "B82";"G82";
	\ar@{-} "B82";"D82";
	\ar@{-} "B82";"E82";
	\ar@{} "B82";"F82";
	\ar@{} "C82";"E82";
	\ar@{} "C82";"F82";
	\ar@{} "C82";"G82";
	\ar@{} "D82";"F82";
	\ar@{} "D82";"G82";
	\ar@{} "E82";"G82";	

	\ar@{} (0, 0);(42, -162) *\cir<2pt>{} = "B83";
	\ar@{-} "B83";(40, -166) *\cir<2pt>{} = "C83";
	\ar@{-} "C83";(42, -170) *\cir<2pt>{} = "D83";
	\ar@{-} "D83";(46, -170) *\cir<2pt>{} = "E83";
	\ar@{-} "E83";(48, -166) *\cir<2pt>{} = "F83";
	\ar@{} "F83";(46, -162) *\cir<2pt>{} = "G83";
	\ar@{-} "B83";"G83";
	\ar@{-} "B83";"D83";
	\ar@{} "B83";"E83";
	\ar@{} "B83";"F83";
	\ar@{-} "C83";"E83";
	\ar@{} "C83";"F83";
	\ar@{} "C83";"G83";
	\ar@{} "D83";"F83";
	\ar@{} "D83";"G83";
	\ar@{-} "E83";"G83";	

	\ar@{} (0, 0);(54, -162) *\cir<2pt>{} = "B84";
	\ar@{-} "B84";(52, -166) *\cir<2pt>{} = "C84";
	\ar@{-} "C84";(54, -170) *\cir<2pt>{} = "D84";
	\ar@{-} "D84";(58, -170) *\cir<2pt>{} = "E84";
	\ar@{} "E84";(60, -166) *\cir<2pt>{} = "F84";
	\ar@{-} "F84";(58, -162) *\cir<2pt>{} = "G84";
	\ar@{-} "B84";"G84";
	\ar@{} "B84";"D84";
	\ar@{-} "B84";"E84";
	\ar@{} "B84";"F84";
	\ar@{-} "C84";"E84";
	\ar@{} "C84";"F84";
	\ar@{} "C84";"G84";
	\ar@{} "D84";"F84";
	\ar@{-} "D84";"G84";
	\ar@{} "E84";"G84";	

	\ar@{} (0, 0);(66, -162) *\cir<2pt>{} = "B85";
	\ar@{-} "B85";(64, -166) *\cir<2pt>{} = "C85";
	\ar@{-} "C85";(66, -170) *\cir<2pt>{} = "D85";
	\ar@{-} "D85";(70, -170) *\cir<2pt>{} = "E85";
	\ar@{-} "E85";(72, -166) *\cir<2pt>{} = "F85";
	\ar@{-} "F85";(70, -162) *\cir<2pt>{} = "G85";
	\ar@{} "B85";"G85";
	\ar@{} "B85";"D85";
	\ar@{} "B85";"E85";
	\ar@{-} "B85";"F85";
	\ar@{} "C85";"E85";
	\ar@{} "C85";"F85";
	\ar@{-} "C85";"G85";
	\ar@{-} "D85";"F85";
	\ar@{} "D85";"G85";
	\ar@{} "E85";"G85";	

	\ar@{} (0, 0);(78, -162) *\cir<2pt>{} = "B86";
	\ar@{-} "B86";(76, -166) *\cir<2pt>{} = "C86";
	\ar@{-} "C86";(78, -170) *\cir<2pt>{} = "D86";
	\ar@{-} "D86";(82, -170) *\cir<2pt>{} = "E86";
	\ar@{-} "E86";(84, -166) *\cir<2pt>{} = "F86";
	\ar@{-} "F86";(82, -162) *\cir<2pt>{} = "G86";
	\ar@{-} "B86";"G86";
	\ar@{-} "B86";"D86";
	\ar@{} "B86";"E86";
	\ar@{} "B86";"F86";
	\ar@{} "C86";"E86";
	\ar@{} "C86";"F86";
	\ar@{} "C86";"G86";
	\ar@{} "D86";"F86";
	\ar@{} "D86";"G86";
	\ar@{-} "E86";"G86";	

	\ar@{} (0, 0);(90, -162) *\cir<2pt>{} = "B87";
	\ar@{-} "B87";(88, -166) *\cir<2pt>{} = "C87";
	\ar@{-} "C87";(90, -170) *\cir<2pt>{} = "D87";
	\ar@{-} "D87";(94, -170) *\cir<2pt>{} = "E87";
	\ar@{-} "E87";(96, -166) *\cir<2pt>{} = "F87";
	\ar@{} "F87";(94, -162) *\cir<2pt>{} = "G87";
	\ar@{-} "B87";"G87";
	\ar@{} "B87";"D87";
	\ar@{-} "B87";"E87";
	\ar@{} "B87";"F87";
	\ar@{} "C87";"E87";
	\ar@{} "C87";"F87";
	\ar@{} "C87";"G87";
	\ar@{} "D87";"F87";
	\ar@{-} "D87";"G87";
	\ar@{-} "E87";"G87";	

	\ar@{} (0, 0);(102, -162) *\cir<2pt>{} = "B88";
	\ar@{-} "B88";(100, -166) *\cir<2pt>{} = "C88";
	\ar@{-} "C88";(102, -170) *\cir<2pt>{} = "D88";
	\ar@{-} "D88";(106, -170) *\cir<2pt>{} = "E88";
	\ar@{-} "E88";(108, -166) *\cir<2pt>{} = "F88";
	\ar@{-} "F88";(106, -162) *\cir<2pt>{} = "G88";
	\ar@{} "B88";"G88";
	\ar@{} "B88";"D88";
	\ar@{} "B88";"E88";
	\ar@{-} "B88";"F88";
	\ar@{} "C88";"E88";
	\ar@{-} "C88";"F88";
	\ar@{-} "C88";"G88";
	\ar@{} "D88";"F88";
	\ar@{} "D88";"G88";
	\ar@{} "E88";"G88";	

	\ar@{} (0, 0);(114, -162) *\cir<2pt>{} = "B89";
	\ar@{-} "B89";(112, -166) *\cir<2pt>{} = "C89";
	\ar@{-} "C89";(114, -170) *\cir<2pt>{} = "D89";
	\ar@{-} "D89";(118, -170) *\cir<2pt>{} = "E89";
	\ar@{-} "E89";(120, -166) *\cir<2pt>{} = "F89";
	\ar@{-} "F89";(118, -162) *\cir<2pt>{} = "G89";
	\ar@{-} "B89";"G89";
	\ar@{} "B89";"D89";
	\ar@{-} "B89";"E89";
	\ar@{} "B89";"F89";
	\ar@{} "C89";"E89";
	\ar@{} "C89";"F89";
	\ar@{} "C89";"G89";
	\ar@{} "D89";"F89";
	\ar@{-} "D89";"G89";
	\ar@{} "E89";"G89";	

	\ar@{} (0, 0);(126, -162) *\cir<2pt>{} = "B90";
	\ar@{-} "B90";(124, -166) *\cir<2pt>{} = "C90";
	\ar@{-} "C90";(126, -170) *\cir<2pt>{} = "D90";
	\ar@{-} "D90";(130, -170) *\cir<2pt>{} = "E90";
	\ar@{-} "E90";(132, -166) *\cir<2pt>{} = "F90";
	\ar@{-} "F90";(130, -162) *\cir<2pt>{} = "G90";
	\ar@{-} "B90";"G90";
	\ar@{-} "B90";"D90";
	\ar@{} "B90";"E90";
	\ar@{} "B90";"F90";
	\ar@{} "C90";"E90";
	\ar@{-} "C90";"F90";
	\ar@{} "C90";"G90";
	\ar@{} "D90";"F90";
	\ar@{} "D90";"G90";
	\ar@{-} "E90";"G90";	

	\ar@{} (0, 0);(18, -174) *\cir<2pt>{} = "B91";
	\ar@{-} "B91";(16, -178) *\cir<2pt>{} = "C91";
	\ar@{-} "C91";(18, -182) *\cir<2pt>{} = "D91";
	\ar@{-} "D91";(22, -182) *\cir<2pt>{} = "E91";
	\ar@{-} "E91";(24, -178) *\cir<2pt>{} = "F91";
	\ar@{-} "F91";(22, -174) *\cir<2pt>{} = "G91";
	\ar@{-} "B91";"G91";
	\ar@{-} "B91";"D91";
	\ar@{} "B91";"E91";
	\ar@{} "B91";"F91";
	\ar@{} "C91";"E91";
	\ar@{-} "C91";"F91";
	\ar@{} "C91";"G91";
	\ar@{-} "D91";"F91";
	\ar@{} "D91";"G91";
	\ar@{} "E91";"G91";	

	\ar@{} (0, 0);(30, -174) *\cir<2pt>{} = "B92";
	\ar@{-} "B92";(28, -178) *\cir<2pt>{} = "C92";
	\ar@{-} "C92";(30, -182) *\cir<2pt>{} = "D92";
	\ar@{-} "D92";(34, -182) *\cir<2pt>{} = "E92";
	\ar@{-} "E92";(36, -178) *\cir<2pt>{} = "F92";
	\ar@{-} "F92";(34, -174) *\cir<2pt>{} = "G92";
	\ar@{} "B92";"G92";
	\ar@{} "B92";"D92";
	\ar@{} "B92";"E92";
	\ar@{-} "B92";"F92";
	\ar@{-} "C92";"E92";
	\ar@{} "C92";"F92";
	\ar@{-} "C92";"G92";
	\ar@{-} "D92";"F92";
	\ar@{} "D92";"G92";
	\ar@{} "E92";"G92";	

	\ar@{} (0, 0);(42, -174) *\cir<2pt>{} = "B93";
	\ar@{-} "B93";(40, -178) *\cir<2pt>{} = "C93";
	\ar@{-} "C93";(42, -182) *\cir<2pt>{} = "D93";
	\ar@{-} "D93";(46, -182) *\cir<2pt>{} = "E93";
	\ar@{-} "E93";(48, -178) *\cir<2pt>{} = "F93";
	\ar@{} "F93";(46, -174) *\cir<2pt>{} = "G93";
	\ar@{-} "B93";"G93";
	\ar@{-} "B93";"D93";
	\ar@{} "B93";"E93";
	\ar@{} "B93";"F93";
	\ar@{-} "C93";"E93";
	\ar@{} "C93";"F93";
	\ar@{} "C93";"G93";
	\ar@{} "D93";"F93";
	\ar@{-} "D93";"G93";
	\ar@{-} "E93";"G93";	

	\ar@{} (0, 0);(54, -174) *\cir<2pt>{} = "B94";
	\ar@{-} "B94";(52, -178) *\cir<2pt>{} = "C94";
	\ar@{-} "C94";(54, -182) *\cir<2pt>{} = "D94";
	\ar@{-} "D94";(58, -182) *\cir<2pt>{} = "E94";
	\ar@{-} "E94";(60, -178) *\cir<2pt>{} = "F94";
	\ar@{-} "F94";(58, -174) *\cir<2pt>{} = "G94";
	\ar@{-} "B94";"G94";
	\ar@{} "B94";"D94";
	\ar@{-} "B94";"E94";
	\ar@{} "B94";"F94";
	\ar@{} "C94";"E94";
	\ar@{} "C94";"F94";
	\ar@{} "C94";"G94";
	\ar@{} "D94";"F94";
	\ar@{-} "D94";"G94";
	\ar@{-} "E94";"G94";	

	\ar@{} (0, 0);(66, -174) *\cir<2pt>{} = "B95";
	\ar@{-} "B95";(64, -178) *\cir<2pt>{} = "C95";
	\ar@{-} "C95";(66, -182) *\cir<2pt>{} = "D95";
	\ar@{-} "D95";(70, -182) *\cir<2pt>{} = "E95";
	\ar@{-} "E95";(72, -178) *\cir<2pt>{} = "F95";
	\ar@{} "F95";(70, -174) *\cir<2pt>{} = "G95";
	\ar@{-} "B95";"G95";
	\ar@{} "B95";"D95";
	\ar@{-} "B95";"E95";
	\ar@{} "B95";"F95";
	\ar@{-} "C95";"E95";
	\ar@{} "C95";"F95";
	\ar@{} "C95";"G95";
	\ar@{} "D95";"F95";
	\ar@{-} "D95";"G95";
	\ar@{-} "E95";"G95";	

	\ar@{} (0, 0);(78, -174) *\cir<2pt>{} = "B96";
	\ar@{-} "B96";(76, -178) *\cir<2pt>{} = "C96";
	\ar@{-} "C96";(78, -182) *\cir<2pt>{} = "D96";
	\ar@{-} "D96";(82, -182) *\cir<2pt>{} = "E96";
	\ar@{-} "E96";(84, -178) *\cir<2pt>{} = "F96";
	\ar@{-} "F96";(82, -174) *\cir<2pt>{} = "G96";
	\ar@{-} "B96";"G96";
	\ar@{-} "B96";"D96";
	\ar@{} "B96";"E96";
	\ar@{} "B96";"F96";
	\ar@{} "C96";"E96";
	\ar@{} "C96";"F96";
	\ar@{-} "C96";"G96";
	\ar@{} "D96";"F96";
	\ar@{} "D96";"G96";
	\ar@{-} "E96";"G96";	

	\ar@{} (0, 0);(90, -174) *\cir<2pt>{} = "B97";
	\ar@{-} "B97";(88, -178) *\cir<2pt>{} = "C97";
	\ar@{-} "C97";(90, -182) *\cir<2pt>{} = "D97";
	\ar@{-} "D97";(94, -182) *\cir<2pt>{} = "E97";
	\ar@{-} "E97";(96, -178) *\cir<2pt>{} = "F97";
	\ar@{-} "F97";(94, -174) *\cir<2pt>{} = "G97";
	\ar@{-} "B97";"G97";
	\ar@{-} "B97";"D97";
	\ar@{-} "B97";"E97";
	\ar@{} "B97";"F97";
	\ar@{} "C97";"E97";
	\ar@{-} "C97";"F97";
	\ar@{} "C97";"G97";
	\ar@{} "D97";"F97";
	\ar@{} "D97";"G97";
	\ar@{} "E97";"G97";	

	\ar@{} (0, 0);(102, -174) *\cir<2pt>{} = "B98";
	\ar@{-} "B98";(100, -178) *\cir<2pt>{} = "C98";
	\ar@{-} "C98";(102, -182) *\cir<2pt>{} = "D98";
	\ar@{-} "D98";(106, -182) *\cir<2pt>{} = "E98";
	\ar@{-} "E98";(108, -178) *\cir<2pt>{} = "F98";
	\ar@{-} "F98";(106, -174) *\cir<2pt>{} = "G98";
	\ar@{-} "B98";"G98";
	\ar@{-} "B98";"D98";
	\ar@{} "B98";"E98";
	\ar@{} "B98";"F98";
	\ar@{} "C98";"E98";
	\ar@{} "C98";"F98";
	\ar@{-} "C98";"G98";
	\ar@{} "D98";"F98";
	\ar@{-} "D98";"G98";
	\ar@{} "E98";"G98";	

	\ar@{} (0, 0);(114, -174) *\cir<2pt>{} = "B99";
	\ar@{-} "B99";(112, -178) *\cir<2pt>{} = "C99";
	\ar@{-} "C99";(114, -182) *\cir<2pt>{} = "D99";
	\ar@{-} "D99";(118, -182) *\cir<2pt>{} = "E99";
	\ar@{-} "E99";(120, -178) *\cir<2pt>{} = "F99";
	\ar@{-} "F99";(118, -174) *\cir<2pt>{} = "G99";
	\ar@{-} "B99";"G99";
	\ar@{} "B99";"D99";
	\ar@{} "B99";"E99";
	\ar@{-} "B99";"F99";
	\ar@{-} "C99";"E99";
	\ar@{} "C99";"F99";
	\ar@{-} "C99";"G99";
	\ar@{} "D99";"F99";
	\ar@{} "D99";"G99";
	\ar@{-} "E99";"G99";	

	\ar@{} (0, 0);(126, -174) *\cir<2pt>{} = "B100";
	\ar@{-} "B100";(124, -178) *\cir<2pt>{} = "C100";
	\ar@{-} "C100";(126, -182) *\cir<2pt>{} = "D100";
	\ar@{-} "D100";(130, -182) *\cir<2pt>{} = "E100";
	\ar@{-} "E100";(132, -178) *\cir<2pt>{} = "F100";
	\ar@{-} "F100";(130, -174) *\cir<2pt>{} = "G100";
	\ar@{-} "B100";"G100";
	\ar@{} "B100";"D100";
	\ar@{} "B100";"E100";
	\ar@{-} "B100";"F100";
	\ar@{-} "C100";"E100";
	\ar@{} "C100";"F100";
	\ar@{-} "C100";"G100";
	\ar@{-} "D100";"F100";
	\ar@{} "D100";"G100";
	\ar@{} "E100";"G100";	

	\ar@{} (0, 0);(18, -186) *\cir<2pt>{} = "B101";
	\ar@{-} "B101";(16, -190) *\cir<2pt>{} = "C101";
	\ar@{-} "C101";(18, -194) *\cir<2pt>{} = "D101";
	\ar@{-} "D101";(22, -194) *\cir<2pt>{} = "E101";
	\ar@{-} "E101";(24, -190) *\cir<2pt>{} = "F101";
	\ar@{-} "F101";(22, -186) *\cir<2pt>{} = "G101";
	\ar@{-} "B101";"G101";
	\ar@{} "B101";"D101";
	\ar@{-} "B101";"E101";
	\ar@{} "B101";"F101";
	\ar@{-} "C101";"E101";
	\ar@{} "C101";"F101";
	\ar@{} "C101";"G101";
	\ar@{} "D101";"F101";
	\ar@{-} "D101";"G101";
	\ar@{-} "E101";"G101";	

	\ar@{} (0, 0);(30, -186) *\cir<2pt>{} = "B102";
	\ar@{-} "B102";(28, -190) *\cir<2pt>{} = "C102";
	\ar@{-} "C102";(30, -194) *\cir<2pt>{} = "D102";
	\ar@{-} "D102";(34, -194) *\cir<2pt>{} = "E102";
	\ar@{-} "E102";(36, -190) *\cir<2pt>{} = "F102";
	\ar@{-} "F102";(34, -186) *\cir<2pt>{} = "G102";
	\ar@{-} "B102";"G102";
	\ar@{} "B102";"D102";
	\ar@{} "B102";"E102";
	\ar@{-} "B102";"F102";
	\ar@{-} "C102";"E102";
	\ar@{} "C102";"F102";
	\ar@{} "C102";"G102";
	\ar@{-} "D102";"F102";
	\ar@{} "D102";"G102";
	\ar@{-} "E102";"G102";	

	\ar@{} (0, 0);(42, -186) *\cir<2pt>{} = "B103";
	\ar@{-} "B103";(40, -190) *\cir<2pt>{} = "C103";
	\ar@{-} "C103";(42, -194) *\cir<2pt>{} = "D103";
	\ar@{-} "D103";(46, -194) *\cir<2pt>{} = "E103";
	\ar@{-} "E103";(48, -190) *\cir<2pt>{} = "F103";
	\ar@{-} "F103";(46, -186) *\cir<2pt>{} = "G103";
	\ar@{-} "B103";"G103";
	\ar@{} "B103";"D103";
	\ar@{-} "B103";"E103";
	\ar@{} "B103";"F103";
	\ar@{} "C103";"E103";
	\ar@{-} "C103";"F103";
	\ar@{-} "C103";"G103";
	\ar@{} "D103";"F103";
	\ar@{} "D103";"G103";
	\ar@{-} "E103";"G103";	

	\ar@{} (0, 0);(54, -186) *\cir<2pt>{} = "B104";
	\ar@{-} "B104";(52, -190) *\cir<2pt>{} = "C104";
	\ar@{-} "C104";(54, -194) *\cir<2pt>{} = "D104";
	\ar@{-} "D104";(58, -194) *\cir<2pt>{} = "E104";
	\ar@{-} "E104";(60, -190) *\cir<2pt>{} = "F104";
	\ar@{-} "F104";(58, -186) *\cir<2pt>{} = "G104";
	\ar@{-} "B104";"G104";
	\ar@{} "B104";"D104";
	\ar@{-} "B104";"E104";
	\ar@{} "B104";"F104";
	\ar@{-} "C104";"E104";
	\ar@{-} "C104";"F104";
	\ar@{} "C104";"G104";
	\ar@{-} "D104";"F104";
	\ar@{} "D104";"G104";
	\ar@{} "E104";"G104";	

	\ar@{} (0, 0);(66, -186) *\cir<2pt>{} = "B105";
	\ar@{-} "B105";(64, -190) *\cir<2pt>{} = "C105";
	\ar@{-} "C105";(66, -194) *\cir<2pt>{} = "D105";
	\ar@{-} "D105";(70, -194) *\cir<2pt>{} = "E105";
	\ar@{-} "E105";(72, -190) *\cir<2pt>{} = "F105";
	\ar@{-} "F105";(70, -186) *\cir<2pt>{} = "G105";
	\ar@{-} "B105";"G105";
	\ar@{} "B105";"D105";
	\ar@{-} "B105";"E105";
	\ar@{} "B105";"F105";
	\ar@{} "C105";"E105";
	\ar@{-} "C105";"F105";
	\ar@{} "C105";"G105";
	\ar@{} "D105";"F105";
	\ar@{-} "D105";"G105";
	\ar@{-} "E105";"G105";	

	\ar@{} (0, 0);(78, -186) *\cir<2pt>{} = "B106";
	\ar@{-} "B106";(76, -190) *\cir<2pt>{} = "C106";
	\ar@{-} "C106";(78, -194) *\cir<2pt>{} = "D106";
	\ar@{-} "D106";(82, -194) *\cir<2pt>{} = "E106";
	\ar@{-} "E106";(84, -190) *\cir<2pt>{} = "F106";
	\ar@{-} "F106";(82, -186) *\cir<2pt>{} = "G106";
	\ar@{-} "B106";"G106";
	\ar@{} "B106";"D106";
	\ar@{} "B106";"E106";
	\ar@{-} "B106";"F106";
	\ar@{-} "C106";"E106";
	\ar@{} "C106";"F106";
	\ar@{-} "C106";"G106";
	\ar@{-} "D106";"F106";
	\ar@{-} "D106";"G106";
	\ar@{} "E106";"G106";	

	\ar@{} (0, 0);(90, -186) *\cir<2pt>{} = "B107";
	\ar@{-} "B107";(88, -190) *\cir<2pt>{} = "C107";
	\ar@{-} "C107";(90, -194) *\cir<2pt>{} = "D107";
	\ar@{-} "D107";(94, -194) *\cir<2pt>{} = "E107";
	\ar@{-} "E107";(96, -190) *\cir<2pt>{} = "F107";
	\ar@{-} "F107";(94, -186) *\cir<2pt>{} = "G107";
	\ar@{-} "B107";"G107";
	\ar@{-} "B107";"D107";
	\ar@{-} "B107";"E107";
	\ar@{} "B107";"F107";
	\ar@{} "C107";"E107";
	\ar@{-} "C107";"F107";
	\ar@{} "C107";"G107";
	\ar@{} "D107";"F107";
	\ar@{-} "D107";"G107";
	\ar@{-} "E107";"G107";	

	\ar@{} (0, 0);(102, -186) *\cir<2pt>{} = "B108";
	\ar@{-} "B108";(100, -190) *\cir<2pt>{} = "C108";
	\ar@{-} "C108";(102, -194) *\cir<2pt>{} = "D108";
	\ar@{-} "D108";(106, -194) *\cir<2pt>{} = "E108";
	\ar@{-} "E108";(108, -190) *\cir<2pt>{} = "F108";
	\ar@{-} "F108";(106, -186) *\cir<2pt>{} = "G108";
	\ar@{} "B108";"G108";
	\ar@{-} "B108";"D108";
	\ar@{-} "B108";"E108";
	\ar@{-} "B108";"F108";
	\ar@{-} "C108";"E108";
	\ar@{} "C108";"F108";
	\ar@{-} "C108";"G108";
	\ar@{-} "D108";"F108";
	\ar@{} "D108";"G108";
	\ar@{} "E108";"G108";	

	\ar@{} (0, 0);(114, -186) *\cir<2pt>{} = "B109";
	\ar@{-} "B109";(112, -190) *\cir<2pt>{} = "C109";
	\ar@{-} "C109";(114, -194) *\cir<2pt>{} = "D109";
	\ar@{-} "D109";(118, -194) *\cir<2pt>{} = "E109";
	\ar@{-} "E109";(120, -190) *\cir<2pt>{} = "F109";
	\ar@{-} "F109";(118, -186) *\cir<2pt>{} = "G109";
	\ar@{-} "B109";"G109";
	\ar@{} "B109";"D109";
	\ar@{} "B109";"E109";
	\ar@{-} "B109";"F109";
	\ar@{-} "C109";"E109";
	\ar@{-} "C109";"F109";
	\ar@{} "C109";"G109";
	\ar@{-} "D109";"F109";
	\ar@{-} "D109";"G109";
	\ar@{} "E109";"G109";	

	\ar@{} (0, 0);(126, -186) *\cir<2pt>{} = "B110";
	\ar@{-} "B110";(124, -190) *\cir<2pt>{} = "C110";
	\ar@{-} "C110";(126, -194) *\cir<2pt>{} = "D110";
	\ar@{-} "D110";(130, -194) *\cir<2pt>{} = "E110";
	\ar@{-} "E110";(132, -190) *\cir<2pt>{} = "F110";
	\ar@{-} "F110";(130, -186) *\cir<2pt>{} = "G110";
	\ar@{-} "B110";"G110";
	\ar@{} "B110";"D110";
	\ar@{-} "B110";"E110";
	\ar@{} "B110";"F110";
	\ar@{-} "C110";"E110";
	\ar@{-} "C110";"F110";
	\ar@{} "C110";"G110";
	\ar@{-} "D110";"F110";
	\ar@{-} "D110";"G110";
	\ar@{-} "E110";"G110";

	\textcolor{green}{\ar@{-} (-6, -42);(15, -42);}
	\textcolor{green}{\ar@{} (0, 0);(23, -42) *\txt{Chordal};}
	\textcolor{green}{\ar@{-} (30, -42);(131, -42);}
	\textcolor{green}{\ar@{-} (-10, -42);(-10, -142);}
	\textcolor{green}{\ar@{-} (128, -42);(128, -142);}
	\textcolor{green}{\ar@{-} (-13, -142);(127, -142);}

	\textcolor{red}{\ar@{-} (70, -46);(102, -46);}
	\textcolor{red}{\ar@{} (0, 0);(108, -46) *\txt{Closed};}
	\textcolor{red}{\ar@{-} (114, -46);(120, -46);}
	\textcolor{red}{\ar@{-} (66, -46);(66, -138);}
	\textcolor{red}{\ar@{-} (117, -46);(117, -138);}
	\textcolor{red}{\ar@{-} (63, -138);(116, -138);}

\end{xy}

\end{document}